%% file: monadicity_measurability.tex
\documentclass[reqno,twoside]{amsart}
\usepackage{monadicity-measurability}

\input{metadata.tex}

\begin{document}

\input{sections/00_abstract.tex}
\maketitle
\tableofcontents

\input{sections/01_introduction.tex}
\input{sections/02_preliminaries.tex}

\input{sections/03_monadicity.tex}

\input{sections/04_giry_comonads.tex}

\input{sections/05_local_presentability.tex}

\input{sections/06_tensor_product.tex}

\input{sections/07_sigma_cp_instruments.tex}

\bibliographystyle{amsalpha}
\bibliography{references}
\end{document}

%% file: metadata.tex
\title{Monadicity of noncommutative measurable spaces over $C^{\ast}$-algebras}
\paperauthor{1}{Fosco}{Loregian}{F.}{Tallin University of Technology}{Tallinn}{Estonia}{fosco.loregian@taltech.ee}
\paperauthor{2}{Antonio}{Lorenzin}{A.}{}{Trento}{Italy}{a.lorenzin.95@gmail.com}

%% file: sections/00_abstract.tex
\begin{abstract}
	In a recent paper, Tobias Fritz and the second author identified  $\sigma C^*$-algebras as a tentative noncommutative analogue of measurable spaces. 
	We advance this analogy by establishing monadicity of $\sigma C^*$-algebras over $C^*$-algebras and proving several results of interest to both the category theory and operator algebra communities.
	In detail, we show that:
	\begin{itemize}
		\item The associated probabilistic morphisms, $\sigma$-cpu maps, admit a description via a comonad closely related to the Giry monad.
		      Equivalently, $\sigma$-cpu maps can be characterized by a universal property.

		\item The category of $\sigma C^*$-algebras and $\sigma$-homomorphisms is locally $\aleph_1$-presentable.

		\item $\sigma C^*$-algebras have a tensor product induced by the maximal tensor product of $C^*$-algebras and satisfying an analogous universal property.

		\item The notion of CP-instruments, defined as $\sigma$-additive assignments from a $\sigma$-algebra to a space of completely positive maps, can be adapted to this setting by taking as target the set of $\sigma$-normal completely positive maps between $\sigma C^*$-algebras. This yields a bijective correspondence with $\sigma$-cpu maps out of the tensor product of the previous item.
	\end{itemize}
\end{abstract}

%% file: sections/01_introduction.tex
%\anto{to be sent to: Journal of the European Mathematical Society?}

\section{Introduction}
Gelfand duality \cite[\S IV.4]{johnstone82stone} is a famous result in operator algebra establishing that the category of (unital) commutative $C^*$-algebras with $*$-homomorphisms is contravariantly equivalent to the category of compact Hausdorff spaces with continuous maps.
Together with the Riesz--Markov--Kakutani representation theorem
\cite[Theorem 2.14]{rudin1987analysis} --- which ensures that regular Borel probability measures on a compact Hausdorff space $X$ bijectively correspond to states on $C(X)$, i.e., completely positive unital (\cpu) functionals --- this perspective has motivated the adoption of $C^*$-algebras as a natural framework for noncommutative probability.
In recent years, Furber and Jacobs~\cite{furber_jacobs_gelfand} have further supported this adoption by establishing a \emph{probabilistic} Gelfand duality: \cpu maps between (unital) commutative $C^*$-algebras correspond to \emph{continuous} Markov kernels between compact Hausdorff spaces. 
Crucially, however, although disintegration of measures (or, in more naive terms, ``generalized Bayes' rule'') yield Markov kernels, it rarely produces continuous ones.

To overcome this limitation, one generally resorts to von Neumann algebras, often described as ``noncommutative measurable spaces''.
While this slogan is well motivated, it obscures an important distinction:
commutative von Neumann algebras are ``blind to null sets'' because their elements are \emph{equivalence classes} of random variables (modulo almost-everywhere equality).\footnote{A genuine Gelfand--type duality for von Neumann algebras requires a substantial revision of the notion of measurable space; see \cite{pavlov2022duality}.}
Formally, commutative von Neumann algebras correspond to measure algebras, i.e.\ complete Boolean algebras obtained by quotienting $\sigma$-algebras by null sets.

This might initially seem like a harmless technicality, but it has undesirable consequences.
For instance, one may wish to treat the Lebesgue measure on $\bbR$ together with delta measures.
However, singleton sets are null for the Lebesgue measure, so quotienting by null sets erases the information carried by individual points.

An approach tailored to measurable spaces was proposed by Tobias Fritz and the second author in \cite{fritz2025categories}.
In that work, the tentative noncommutative analogue of measurable spaces is given by \emph{\cstars{}s}, i.e.\ monotone $\sigma$-complete \cstar{}s~\cite{saito15monotone}.\footnote{Preliminaries on \cstars{}s are provided in \cref{sec:preliminaries}.
	For now, it suffices to recall that the standard morphisms between \cstars{}s are the \emph{$\sigma$-normal} maps.}
This proposal is supported by a Gelfand--type duality: assigning to every measurable space $X$ the commutative \cstars{} $\Linf(X)$ of bounded measurable functions yields an adjunction.
While this is not an equivalence in full generality, it restricts to one on standard Borel spaces and, beyond second-countability, on Baire measurable spaces~\cite[\S 4.2]{fritz2025categories}.
Moreover, the corresponding probabilistic morphisms, namely $\sigma$-normal \cpu{} maps, are closely related to positive operator-valued measures (POVMs).\footnote{This connection will also reappear in the present work; see \cref{sec:sigmaCP}.}
 
\begin{figure}[ht]
	\centering
	\begin{tikzpicture}[x=3.8cm,y=1.25cm,
			every node/.style={align=center,font=\small,inner sep=4pt},
			>=stealth]
		\node (preliminaries) at (0,0)
			{\cref{sec:preliminaries}\\Preliminaries};
		\node (monadicity) at (0,-1)
			{\cref{sec:monadic}\\Monadicity};
		\node (giry) at (-1,-2)
			{\cref{sec:giry}\\Giry (co)monads};
		\node (accessible) at (0,-2)
			{\cref{sec:accessible}\\Local presentability};
		\node (tensor) at (1,-2)
			{\cref{sec:tensor}\\Tensor product};
		\node (instruments) at (1,-3)
			{\cref{sec:sigmaCP}\\$\sigma$CP-instruments};
		\draw[->] (preliminaries) -- (monadicity);
		\draw[->] (monadicity) -- (giry);
		\draw[->] (monadicity) -- (accessible);
		\draw[->] (monadicity) -- (tensor);
		\draw[->] (tensor) -- (instruments);
	\end{tikzpicture}
	\caption{Section dependencies, with transitive arrows omitted.}\label{fig:secdep}
\end{figure}
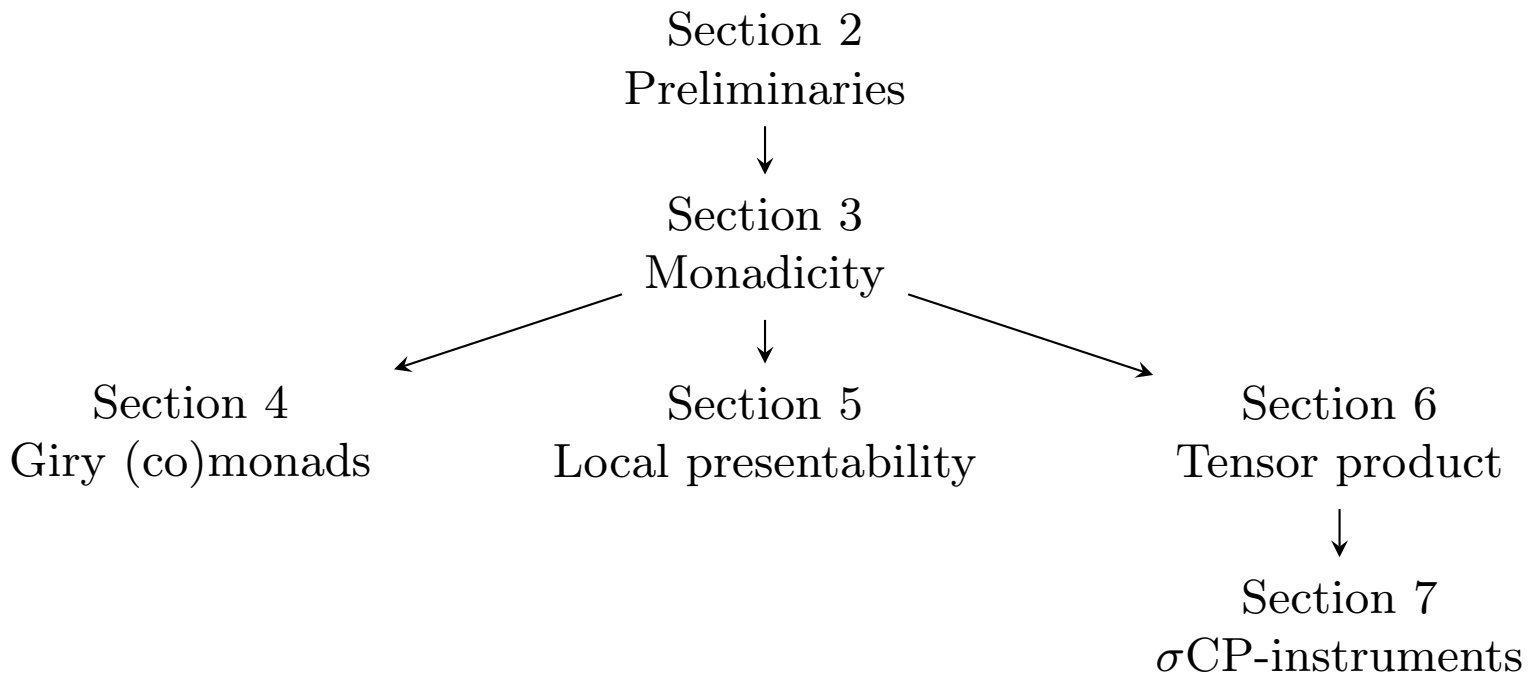

The present work builds on these developments.
The main theorem, which underlies the rest of the paper (see \cref{fig:secdep}), shows that the forgetful functor $U\colon \SigmaC \to \Calg$ from the category of \cstars{}s with $\sigma$-normal $*$-homomorphisms to the category of \cstar{}s with $*$-homomorphisms is \emph{monadic} (\cref{sec:monadic}). 
Roughly speaking, this means that \cstars{}s can be thought of as "algebras over \cstar{}s", and the rest of the paper exploits this description to ensure a number of consequences that are also of interest to readers less familiar with category theory. 
Though we will detail these findings below, we briefly mention that we will be able to describe \scpu{} maps via a universal property and moreover equip \cstars{}s with a tensor product. 
The adopted categorical methods uncover the conceptual unity underlying these results, which would otherwise require separate functional-analytic justifications.

The corresponding monadicity result for von Neumann algebras was established by Furber~\cite{furber2024monad}.
An important distinction is that Furber proves monadicity only for normal $*$-homomorphisms, whereas here we show it holds both for $\sigma$-normal $*$-homomorphisms and for $\sigma$-normal \cpu{} maps (\cref{rem:monadic_cpu}). 
In fact, the present result can be seen as a refinement. 
As observed in \cref{rem:monadicity_W}, Furber's theorem together with \cref{thm:U_monadic} yields that $\Walg \to \SigmaC$, where $\Walg$ denotes the category of von Neumann algebras with normal $*$-homomorphisms, is monadic as well. 
In particular, we obtain a tower of monadic functors
\[  
\Walg \longrightarrow \SigmaC \longrightarrow \Calg
\]
of which the present paper supplies the middle step.

Given the categorical and operator-algebraic nature of this paper, we now describe the results from both perspectives, in order to provide conceptual insight to readers from either community.

\subsection*{From a categorical viewpoint}
\begin{itemize}
	\item Using the fact that $\Calg$ (and its subcategory of commutative \cstar{}s) has a comonad describing \cpu{} maps as coalgebras, we establish the existence of an analogous comonad on $\SigmaC$, describing $\sigma$-normal \cpu{} maps as coalgebras (\cref{sec:giry}).

	      Moreover, the comonad for commutative algebras extends the \emph{Giry monad} of \cite{giry82probability} when restricted to Baire measurable spaces \cite[\S4.2]{fritz2025categories}, which generalize standard Borel spaces beyond second countability.
	\item $\SigmaC$ is locally $\aleph_1$-presentable, and the forgetful functor $U\colon \SigmaC \to \Calg$ is accessible (\cref{sec:accessible}).
	\item $\SigmaC$ is a symmetric monoidal category, and the same holds when extending to $\sigma$-normal \cpu{} maps (\cref{sec:tensor}).
	\item The monoidal category is not closed, but a property similar to monoidal closure holds when one object is commutative (\cref{sec:sigmaCP}).
\end{itemize}
\subsection*{From an operator-algebraic viewpoint}
\begin{itemize}
	\item The \cstar{} of continuous functions over the state space and the universal \cstar{}~\cite[Proposition 8]{kirchberg1998Cstarfromoperator} ensure the existence of \cstars{}s having analogous universal properties (\cref{sec:giry}).

	      Unfolding this, for any \cstars{} $\A$ there exists a unique \cstars{} $\mathcal{P}(\A)$ such that every $\sigma$-normal \cpu{} map $\A \to \B$ factorizes as
	      \[
		      \xymatrix{\A\ar[r]\ar[d] & \B\\ \mathcal{P}(\A)\ar@{-->}[ru]}
	      \]
	      where $\mathcal{P}(\A)\to \B$ is a $\sigma$-normal $*$-homomorphism.
	      Similarly, there is a \emph{commutative} \cstars{} $\mathcal{S}(\A)$ when the universal property is tested only on commutative $\B$.
	\item Roughly speaking, each \cstars{} is built from a small set of `$\aleph_1$-compact' \cstars{}s by taking direct sums and quotients (\cref{sec:accessible}).
	\item \cstars{}s have a tensor product $\kocktensor$ obtained by extending the maximal tensor product $\otimes$ on \cstar{}s (\cref{sec:tensor}).
	      In particular, morphisms of type $\A \kocktensor \B \to \C$ are in bijective correspondence with those of type $\A \otimes \B \to \C$ which are \emph{bi-$\sigma$-normal}, i.e.\ separately $\sigma$-normal in each variable (\cref{def:bisigmanormal}).
	      This resembles the universal property of the binormal tensor product for von Neumann algebras (compare \cite[Proposition 2.1]{okamura2015measurement} with \cref{thm:tensor_shom}), where ``binormal'' means precisely ``normal on each tensor factor''.
	\item The discussion of \cite{okamura2015measurement} regarding \emph{CP-instruments} in the context of von Neumann algebras can be reformulated for \cstars{}s (\cref{sec:sigmaCP}). 
	Below, we offer some insights on this notion and our contributions to this matter.

	CP-instruments are due to Davies and Lewis~\cite{davies1970operational}, who introduced them to treat repeated measurements of \emph{continuous} observables without commutation hypotheses. Ozawa~\cite{ozawa1984quantummeasuring} established their correspondence with measuring processes in the Hilbert space setting. 
	The underlying intuition is that every abstract instrument admits a concrete dynamical realization in much the same way that the GNS construction realizes abstract states on \cstar{}s as concrete vector states in a Hilbert space representation.

	An extension to arbitrary von Neumann algebras was provided in \cite{okamura2015measurement}, where however dynamical realizations are not always ensured. 
	To see where difficulties arise, we mention that this realization requires two steps: first, rewriting CP-instruments using a tensor product; second, verifying whether this tensor product extends to a von Neumann algebra. While the former can be proven in full generality (\cite[Proposition 3.3]{okamura2015measurement}), the latter requires an additional property called the ``normal extension property''.

	In this paper, \cref{thm:sigmaCPinstruments} establishes both of the steps highlighted above by weakening normality (proper to von Neumann algebras) to $\sigma$-normality (proper to \cstars{}s). 
	Crucially, no additional property is required.
\end{itemize}
\begin{paperacknowledgments}
We are grateful to Tallinn University of Technology and Pawel Sobocinski for inviting the second author to visit the institute, thereby sparking the collaboration between the two authors that eventually led to this manuscript.
We thank Tobias Fritz for his comments and insights on a draft of this paper.
We also thank Robert Booth, whose interest inspired the developments of \cref{sec:sigmaCP}.
\begin{paperfunding}
{The first author gratefully acknowledges the support of Tallinn University of Technology and of the Estonian Research Council grant~PRG3215.}
\end{paperfunding}
AI tools (particularly, Claude and Gemini Notebook) were employed to polish writing and ideas; their technical contributions remained marginal and were always verified by the authors, who take full responsibility for any errors in the manuscript.\footnote{Almost certainly, such errors will have a human fingerprint anyway.}
\end{paperacknowledgments}

%% file: sections/02_preliminaries.tex
\section{Prerequisites and preliminaries on \texorpdfstring{\cstars{}s}{σC*-algebras}}\label{sec:preliminaries}
Throughout this paper, the proof techniques rely on machinery from category theory, particularly monads and symmetric monoidal categories (as well as monoidal monads).
We assume that the reader is familiar with basic category theory, including functors, natural transformations, limits and adjunctions. 
Prior exposure to monads is not required, though a general sense of them is helpful; for additional background on monads, readers may consult \cite[Sections~4.1 and~4.2]{borceux94handbook2}.
As for symmetric monoidal categories, the usual intuition of a ``tensor product'' should be sufficient.\footnote{Surveys on monoidal categories are abundant: \cite{selinger} is oriented towards applications, mostly in logic, while \cite{EGNO} is foundational for mathematicians interested in representation theory; \cite{milne_notes,borceux94handbook2,cwm} constitute good introductions that do not dwell too much on the finer details.}

This technology is adopted in the categories of (commutative) \cstar{}s and the subclass of (commutative) \cstars{}s. 

\begin{convention}\label{conv:unital}
Throughout this paper, unless otherwise specified (e.g., in the discussion around \autoref{lem:sigmaC_ideal}), \cstar{}s and \cstars{}s are always unital, and $*$-homomorphisms preserve the unit.
\end{convention} 

While the reader is not required to have a complete background on this theory, some basic knowledge of \cstar{}s may be used without further explanations: this includes, e.g., the notion of positivity and the fact that positive unital maps (and so also $*$-homomorphisms) are contractive, i.e.\ of norm $\le 1$~\cite[Proposition 2.1]{paulsen2002}.

Instead, no prior background on \cstars{}s is assumed.
The reader may, however, consult the books \cite{saito15monotone,pedersen2018automgroups}, as well as the article \cite{fritz2025categories}, where the choice of \cstars{}s is motivated from a probabilistic perspective.

\subsection*{\texorpdfstring{\cstars{}s}{σC*-algebras}}
	Recall that a \cstar{} carries a canonical partial order on its self-adjoint part, in which $a \le b$ means that $b-a$ is of the form $c^*c$.
A \newterm{\cstars{}} (also known as a \emph{monotone $\sigma$-complete \cstar{}}) is a \cstar{} $\A$ where every norm-bounded monotone increasing sequence $(a_n)$, i.e.\ such that $a_n \le a_{n+1}$ and $\norm{a_n}\le \lambda$, has a supremum $\sup_n^{\A} a_n$.\footnote{This requirement holds if and only if we restrict it to sequences of \emph{self-adjoint} elements. 
The only if direction is clear, while the if direction follows by noting that $a\le b$ implies that $a-b$ is self-adjoint, so that $a$ and $b$ have the same imaginary part: therefore, in a monotone increasing sequence all terms share the same imaginary part, so the general case reduces to the self-adjoint one upon subtracting it.
Below, the two conditions are used interchangeably.}
The superscript $\A$ in $\sup_n^{\A}$ is meant for clarity, and may be dropped when the context is sufficiently clear.

The slogan to keep in mind is that passing from \cstar{}s to \cstars{}s is the operator-algebraic counterpart of passing from topological to measurable structure: while commutative \cstar{}s correspond to compact Hausdorff spaces, commutative \cstars{}s correspond to \emph{abstract} measurable spaces, as we shall see at the end of this section. 
The intuition is that we allow a limiting procedure weaker than the uniform norm, analogous to how measurable functions are naturally closed under pointwise limits.

Between \cstars{}s, one generally considers \newterm{$\sigma$-normal} maps, i.e.\ those positive maps that respect this class of suprema.
Explicitly, a positive map $\phi\colon \A \to \B$ is $\sigma$-normal when, for every norm-bounded monotone increasing sequence $(a_n)$ in $\A$, the image sequence $(\phi(a_n))$ is again norm-bounded monotone increasing\footnote{This is immediate for positive maps between {unital} \cstar{}s, see \cite[Proposition 2.1]{paulsen2002}.} and
	\[
		\phi\Bigl(\sup\nolimits_n^{\A} a_n\Bigr) = \sup\nolimits_n^{\B} \phi(a_n).
	\]

The \starhom{}s (respectively, \cpu{}\footnote{This abbreviation stands for \emph{c}ompletely
\emph{p}ositive \emph{u}nital.} maps) that are $\sigma$-normal are called \newterm{\starshom{}s} (respectively, \newterm{\scpu{} maps}).
In this paper, we consider the categories as defined in \cref{tab:categories}.
\begin{table}[ht]
	\begin{tabular}{ccc}
		\textbf{Categories} & Objects    & Morphisms    \\
		\toprule
		$\Calg$             & \cstar{}s  & \starhom{}   \\\midrule
		$\Calgcpu$          & \cstar{}s  & \cpu{} maps  \\\midrule
		$\SigmaC$           & \cstars{}s & \starshom{}s \\\midrule
		$\SigmaCcpu$        & \cstars{}s & \scpu{} maps\\ \bottomrule
	\end{tabular}
	\caption{The categories considered in this paper (recall \cref{conv:unital}). The reader should imagine a parallel table where we add $\comm$ in front of each category name, in order to restrict to commutative algebras. For example, $\cCalg$ is the category of commutative \cstar{}s with \starhom{}s, and $\cSigmaCcpu$ is the category of commutative \cstars{}s with \scpu{} maps.}\label{tab:categories}
\end{table}

A \newterm{$\sigma$-subspace} in a \cstars{} $\A$ is a $*$-subspace $V$ having all suprema of norm-bounded monotone increasing sequences: $(a_n)\subseteq V$ implies $\sup_n^{\A}a_n \in V$. %\footnote{A $*$-subspace is simply a vector space closed under involution: $v\in V$ iff $v^*\in V$.}
Similarly, \newterm{$\sigma$-ideals} are those norm-closed two-sided ideals which happen to be $\sigma$-subspaces, and the \newterm{$\sigma$-closure} of a $*$-subspace $W$ is the smallest $\sigma$-subspace containing $W$.

For the technical developments of \cref{sec:tensor,sec:sigmaCP}, we also introduce \newterm{order limits}.
Briefly, one says that a sequence $(a_n)$ is \newterm{order convergent} in a \cstars{} $\A$ if there are four norm-bounded monotone increasing sequences of self-adjoint elements $(b_{n,k})_n$ for $k\in\lbrace 1,2,3,4\rbrace$ such that $a_n=(b_{n,1} - b_{n,2})+i(b_{n,3}-b_{n,4})$, where $i$ is the imaginary unit.
In this case, the order limit is given by
\[\textstyle
	\olim_n^{\A} a_n \coloneqq \sup_n^{\A} b_{n,1} - \sup_n^{\A} b_{n,2} + i \left(\sup_n^{\A} b_{n,3} - \sup_n^{\A} b_{n,4}\right).
\]
Several crucial properties of these limits are contained in \cite[Section 2.1]{saito15monotone}.

In particular, they are:
\begin{itemize}
	\item well-defined, so that $\olim$ of a sequence does not depend on the choice of the $b_{n,k}$'s: whenever we consider norm-bounded monotone increasing sequences of self-adjoint elements $(b_{n,k})_n$ and $(b_{n,k}')_n$ for $k=\lbrace 1,2,3,4\rbrace$ such that
	      \[
		      a_n=(b_{n,1} - b_{n,2})+i(b_{n,3}-b_{n,4})=(b_{n,1}' - b_{n,2}')+i(b_{n,3}'-b_{n,4}') ,
	      \]
	      then
	      \begin{multline}
								\textstyle
		      \sup_n^{\A} b_{n,1} - \sup_n^{\A} b_{n,2} + i \left(\sup_n^{\A} b_{n,3} - \sup_n^{\A} b_{n,4}\right) =\\ \sup_n^{\A} b_{n,1}' - \sup_n^{\A} b_{n,2}' + i \left(\sup_n^{\A} b_{n,3}' - \sup_n^{\A} b_{n,4}' \right).
							\end{multline}
	\item linear: $\olim^{\A}_n (a_n + a_n')=\olim^{\A}_n (a_n)+ \olim^{\A}_n (a_n')$ and $\olim^{\A}_n(\lambda a_n) = \lambda \olim^{\A}_n a_n$ for all $\lambda \in \bbC$;
	\item $*$-preserving: $\olim^{\A}_n (a_n^*)=\olim^{\A}_n(a_n)^*$;
	\item compatible with multiplication: $\olim^{\A}_n (z a_n w) = z(\olim^{\A}_n a_n)w$ for all $z,w \in \A$.
\end{itemize}
Note that suprema do not respect these properties: for example, $-\sup^{\A}_n a_n $ becomes an infimum $\inf^{\A}_n (-a_n)$.
More instructively, multiplying two monotone increasing sequences may not yield a monotone increasing sequence because already the product of two positive elements is not necessarily positive.
Instead, 
\[ \olim^{\A}_n(a_n \olim^{\A}_m b_m)= \olim^{\A}_m (\olim^{\A}_n a_n b_m)=\olim^{\A}_n (a_n)\olim^{\A}_m(b_m).\]

It is a simple check that \scpu{} maps are exactly those cpu maps that respect order limits.
In contrast, it is currently unclear whether a $\sigma$-subspace is necessarily closed under order limits (see \cite[Section 3.5]{fritz2025categories}).

\subsection*{\texorpdfstring{\bsigma{}s}{Boolean σ-algebras}}
Another crucial aspect of \cstars{}s is their relation to \emph{abstract} measurable spaces, better known as \bsigma{}s.
Briefly, these are Boolean algebras with countable joins and meets, and for this reason they mimic $\sigma$-algebras.
The interested reader may refer to \cite{sikorski69boolean} for additional details.
Nonetheless, this paper does not employ the technology of Boolean algebras, and references will be provided when needed.
The material below is needed primarily in \cref{sec:giry,sec:sigmaCP} and may be skipped otherwise.

\begin{notation}
	We denote Boolean $\sigma$-algebras with $\Sigma$ and their elements with $U,V,U_n, \dots$, to remind the reader of the case of measurable spaces.
	Nonetheless, we use the Boolean algebra notations for operations: $\vee$ denotes the join, $\wedge$ the meet, $\neg$ the negation, $\bot$ the bottom element and $\top$ the top.
\end{notation}

It was shown in \cite[Theorem 4.1.34]{fritz2025categories} that $\cSigmaC$ is equivalent to the category of Boolean $\sigma$-algebras and $\sigma$-homomorphisms (here meaning Boolean algebra homomorphisms respecting countable joins and meets).
The \cstars{} associated to the \bsigma{} $\Sigma$ is denoted by $\Linfbool(\Sigma)$.
The notation of $\Linfbool(\Sigma)$ (``abstract $\Linf$'') is chosen to suggest that $\Linfbool(\Sigma)$ is an abstraction of the `concrete' \cstars{} associated to a measurable space $(X,\Sigma_X)$, given by
\[
	\Linf(X)\coloneqq \lbrace f\colon X \to \bbC \text{ bounded measurable}\rbrace
\]
and equipped with pointwise operations.
In fact, $\Linfbool(\Sigma_X)\cong \Linf(X)$ (\cite[Lemma 4.1.14]{fritz2025categories}).
More explicitly, the elements of $\Linfbool(\Sigma)$ are encoded by analogy with preimage maps of bounded measurable functions:
\[
	\Linfbool(\Sigma)\coloneqq \lbrace f\colon B(\bbC)\to \Sigma \mid f(U)=\top \text{ for some bounded }U \in B(\bbC) \rbrace,
\]
where $B(\bbC)$ denotes the Borel $\sigma$-algebra of $\bbC$.
From a commutative \cstars{} $\A$, one can show that $\A \cong \Linfbool(\Sigma)$, with $\Sigma$ its set of projections (which turns out to be a \bsigma{}).
The reader may refer to \cite[Section 4.1]{fritz2025categories} for further details.

From a probabilistic perspective, the category $\Stoch$ of measurable spaces and Markov kernels embeds fully faithfully (contravariantly) into  $\cSigmaCcpu$ (see \cite[Proposition 4.3.6]{fritz2025categories}). Accordingly, $\sigma$-cpu maps $\Linfbool(\Sigma)\to \Linfbool(\Sigma')$ may be viewed as \emph{abstract Markov kernels}.
While tangential to the main discussion, this perspective is included for readers interested in the connection with classical probability.

%% file: sections/03_monadicity.tex
\section{Monadicity of \texorpdfstring{\cstars{}s}{σC*-algebras} over \texorpdfstring{\cstar{}s}{C*-algebras}}\label{sec:monadic}
	This section proves the central result of the paper: the forgetful functor $$U \colon \SigmaC \to \Calg$$ is monadic (\cref{thm:U_monadic}), and the associated monad constructs the \emph{Pedersen--Baire envelope} of a \cstar.
	We then record the two consequences used repeatedly below: $\SigmaC$ is complete with limits created by $U$ (\cref{cor:complete}), and $T$ preserves reflexive coequalizers, whence $\SigmaC$ is cocomplete (\cref{prop:cocomplete}).

We recall the following adjunction, first shown in \cite{wright76regular}. The present formulation is taken from \cite[Equation (3.3.4)]{fritz2025categories}.
\begin{theorem}\label{thm:adj_sigma}
	There is an adjunction
	\begin{equation*}
		% \begin{tikzcd}
		%  \Calg \ar[rr,phantom, "\bot"]\ar[rr,bend left=20,"F"] && \SigmaC \ar[ll,,"U",bend left=20]
		% \end{tikzcd}
		\adPair{\Calg}{\SigmaC}{F}{U}
	\end{equation*}
	where $U$ forgets $\sigma$-normality, and $F$ sends a \cstar{} $\A$ to its \emph{Pedersen--Baire envelope}.
\end{theorem}
Explicitly, the \newterm{Pedersen--Baire envelope} $\A^{\infty}$ of a \cstar{} $\A$ is the $\sigma$-closure of $\A$ inside its enveloping \wstar{} $\A^{**}$~\cite[Definition 3.3.1]{fritz2025categories} or, equivalently, the $\sigma$-closure of $\A$ in $B(H_u)$ for the universal representation $(H_u,\pi_u)$~\cite[Definition 2.2.23]{saito15monotone}. Diagrammatically,
\[
	\xymatrix{
		\A \, \ar@{^{(}->}[r] & \A^{\infty}\,\ar@{^{(}->}[r]^{(\star)}& \A^{**}},
\]
where the right-hand side inclusion $(\star)$ is $\sigma$-normal.
The content of the adjunction of \cref{thm:adj_sigma} is exactly the universal property of $\A^{\infty}$: every $*$-homomorphism from $\A$ into a \cstars{} extends uniquely, and $\sigma$-normally, along $\A \hookrightarrow \A^{\infty}$. While the notation $\A^{\infty}$ is already used in the literature, in this paper we instead opt for $UF\A$ and $F\A$, mindful of the associated adjunction.

\begin{notation}
	Below, $F$ and $U$ always mean the two functors of the previous statement; the associated monad $UF$ will be denoted by $T$ for simplicity; $\epsilon \colon FU \to \id_{\SigmaC}$ and $\eta\colon \id_{\Calg} \to UF$ will denote the counit and unit, respectively, and $\mu\coloneqq U\epsilon F$ is the monad multiplication.
\end{notation}

We now wish to show that $U$ is \emph{monadic}. To say that $U$ is monadic is to say that a \cstars{} is \emph{nothing more than} a \cstar{} equipped with a coherent choice of the countable suprema forced by the Pedersen--Baire envelope: the additional structure is algebraic over $\Calg$, carried by operations and equations rather than by a property. This is what buys, in the sections that follow, limits created by $U$, free objects, a tensor product of algebras (\cref{sec:tensor}) and a presentability theorem (\cref{sec:accessible}), none of which is visible from the definition of a \cstars{}.

	A warning on terminology may spare the reader a search: a substantial part of the literature on monadicity, including \cite{barr2005toposes} and \cite{linton1969coequalizers} as cited below, says \emph{triple} and \emph{tripleable} where we say monad and monadic; \cite{barr2005toposes} does not use the word ``monad'' at all.
	
Monadicity of $U$ can be argued directly, showing that the unique comparison functor $K\colon \SigmaC \to \EM(T)$ deduced from \cite[p. 222]{macdonald2003aspects} is fully faithful and essentially surjective (the two halves of the proof of \cref{thm:U_monadic} below), thus an equivalence of categories.
 The functor $K$ is the only one filling the commutative diagram,
\[
	\xymatrix{
		\SigmaC \ar[rr]^{K} \ar[dr]_{U} && \EM(T) \ar[dl] \\
		& \Calg
	}
\]
where the unlabeled arrow is a forgetful functor.
Explicitly, $K$ sends a \cstars{} $\A$ to the carrier given by $U\A$ equipped with the algebra structure given by the counit $U\epsilon_{\A} \colon TU\A \to U\A$.
\begin{theorem}\label{thm:U_monadic}
	$U$ is monadic.
\end{theorem}
\begin{proof}
	We divide the proof in two steps.
	\setlist[description]{font=\normalfont\itshape}
	\begin{description}
		\item[{$K$ is fully faithful}]   By a proposition due to Beck (see, e.g., \cite[p.~226]{macdonald2003aspects}), $K$ is fully faithful if and only if $\epsilon_{\B}$ is a regular epimorphism.
		      We seek to show that this is always the case.

		      First of all, we recall that $\epsilon_{\B}\colon FU\B \to \B$ is a surjective map in $\SigmaC$.
		      This was proven in \cite[Theorem 5.4.5]{saito15monotone} (see \cite[Theorem 3.3.6]{fritz2025categories} for the explicit statement).
		      As is standard for algebras, every surjective map in fact identifies its image with the quotient via its kernel, i.e. $\B \cong FU\B / \ker(\epsilon_{\B})$.

		      We are therefore left to show that this allows one to write $\epsilon_{\B}$ as a coequalizer. To this end, we consider two cases. If $\ker(\epsilon_{\B})$ does not have a unit, we consider the unitization $\ker(\epsilon_{\B})_1$ of $\ker(\epsilon_{\B})$, which is a \cstars{} as shown in \cite[Proposition 2.2.10]{saito15monotone}.
		      There are two obvious maps $\ker(\epsilon_{\B})_1 \to FU\B$.
		      The first one, which we call $i$, sends $\lambda 1+x$, where $\lambda \in \bbC$ and $x\in \ker(\epsilon_{\B})$, to itself via the inclusion $\ker(\epsilon_{\B}) \subset FU\B $, since $FU\B$ is unital. The second one, denoted by $\ell$, forgets $x$: in other words, it sends $\lambda 1 +x$ to $\lambda 1$. This is well-defined as $\ker(\epsilon_{\B})$ is a (maximal) ideal of its unitization $\ker(\epsilon_{\B})_1$.
		      It is now a simple check that $\epsilon_{\B}$ is the coequalizer of
		      \[
			      \xymatrix{
				      \ker(\epsilon_{\B})_1
				      \ar@<.7ex>[r]^{i}
				      \ar@<-.7ex>[r]_{\ell}
				      & FU\B
			      }
		      \]
		      since the image of their difference $i-\ell$ exactly corresponds to $\ker(\epsilon_{\B})$.
		      If instead $\ker(\epsilon_{\B})$ has a unit, we consider $\ker(\epsilon_{\B})\oplus \bbC$. In this case, we take $i\colon \ker(\epsilon_{\B})\oplus \bbC \to FU\B$ sending $(x,\lambda)$ to $x+\lambda (1-e)$, where $e$ denotes the unit of $\ker(\epsilon_{\B})$, and $\ell$ the composition $\ker(\epsilon_{\B})\oplus \bbC \to \bbC\to FU\B$, where the first map simply sends $(x,\lambda)$ to $\lambda$.
		\item[{$K$ is essentially surjective}] Let us consider an algebra in $\EM(T)$. By definition, this is given by a \cstar{} $\A$ together with a morphism $\phi\colon T\A\to \A$ compatible with the monad multiplication and unit.
		      Explicitly,
		      \begin{equation}\label{eq:EM_algebra}
			      \phi \eta_{\A}=\id_{\A}\qquad\text{and}\qquad\phi T(\phi)= \phi \mu_{\A},
		      \end{equation}
		      where $\eta_{\A}\colon \A \to T\A$ denotes the standard inclusion, while $\mu_{\A}=U(\epsilon_{F\A})\colon TT\A \to T\A$ is the $\sigma$-normal map such that $\mu_{\A}\eta_{T\A}=\id_{T\A}$.

		      Let us show that $\A$ satisfying \eqref{eq:EM_algebra} is a \cstars{} and moreover that $\phi$ is $\sigma$-normal.
		      In other words, $\A = U\A$ and $\phi=U\phi$.
		      Let $(a_n)$ be a norm-bounded monotone increasing sequence in $\A$.
		      Since $\eta_{\A}$ preserves the positive ordering, the image of the sequence $(\eta_{\A}(a_n))$ is a norm-bounded monotone increasing sequence in $T\A$, which is a \cstars{}.
		      We can then consider $\phi(\sup_n^{T\A}(\eta_{\A}(a_n)))\in \A$.
		      Now, given any $b\ge a_n$ for all $n$, $\eta_{\A}(b) \ge \sup_n^{T\A}(\eta_{\A} (a_n))$ and then
		      \[
			      b=\phi(\eta_{\A}(b)) \ge \textstyle\phi\left(\sup_n^{T\A}(\eta_{\A}(a_n)\right)
		      \]
		      Finally, $\sup_n^{\A} a_n $ exists and is given by $\phi (\sup_n^{T\A}(\eta_{\A} (a_n)))$.
		      To show that $\phi$ is $\sigma$-normal, we use the identity given by the algebra multiplication: consider $\sup^{TT\A}_n \eta_{T\A} (x_n)\in TT\A$ for a norm-bounded monotone increasing sequence $(x_n)\subseteq T\A$.
		      Then
			      \begin{align}
				      \textstyle\phi\left(\sup_n^{T\A} \eta_{\A} \phi(x_n)\right) & = \textstyle\phi T(\phi) \left(\sup_n^{TT\A} \eta_{T\A} (x_n)\right)  \tag{definition of $T(\phi)$}\\
				                                                        & = \textstyle\phi \mu_{\A} \left(\sup_n^{TT\A}\eta_{T\A} (x_n)\right) \tag{\ref{eq:EM_algebra}}\\
				                                                        & = \textstyle\phi \left(\sup_n^{T\A} x_n\right)                       \tag{$\sigma$-normality of $\mu$, and $\mu \circ \eta T =\id_{T\A}$}
			      \end{align}
		      Since we just established that $\sup_n^{\A} a_n = \phi (\sup_n^{T\A} (\eta_{\A}(a_n)))$, the left-hand side is equal to $\sup^{\A}_n (\phi(x_n))$.
		      Finally, $\sup^{\A}_n (\phi(x_n)) = \phi(\sup^{T\A}_n x_n)$, as wanted.

		      We now claim that $\phi$ is in fact $\epsilon_{\A}$, from which essential surjectivity of $K$ follows.
		      To this end, consider the set
		      \[M\coloneqq \lbrace a\in T\A \mid \phi(a)=\epsilon_{\A}(a)\rbrace, \]which certainly contains $\eta_{\A} (\A)$ by the equations $\phi \eta_{\A}=\id_{\A}$ and $U(\epsilon_{\A}) \eta_{\A}=\id_{\A}$.
		      Since both $\phi$ and $\epsilon_{\A}$ respect suprema, $M$ is also closed under suprema.
		      Similarly, one also proves that $M$ is a $*$-subspace, hence it is a $\sigma$-subspace.
		      As the $\sigma$-closure of $\eta_{\A} (\A)$ in $T\A$ is the whole $T\A$ (see the discussion after \cref{thm:adj_sigma}), we obtain that $M=T\A$, thus concluding the proof.\qedhere
	\end{description}
\end{proof}

By \cite[Proposition 4.3.1]{borceux94handbook2}, since $\Calg$ has all limits~\cite[Section 2]{khoshkam2002categorical}, \cref{thm:U_monadic} yields the following.

\begin{corollary}\label{cor:complete}
	$\SigmaC$ is complete, and limits are created, in the sense of \cite[\S V.1]{cwm}, by $U$.
\end{corollary}

\begin{remark}[Monadicity for cpu maps, and commutative algebras]\label{rem:monadic_cpu}
	The adjunction of \cref{thm:adj_sigma} can be extended to \cpu{} maps via the \emph{Stinespring dilation theorem}, as shown in \cite[Corollary 3.3.8]{fritz2025categories}.
	Explicitly, we have
	\[
		\adPair{\Calgcpu}{\SigmaCcpu}{\bar F}{\bar U}
	\]
	extending the one in \cref{thm:adj_sigma}.
	The reasoning of \cref{thm:U_monadic} can now be copied verbatim to ensure that the extended $\bar U $ is monadic as well: in particular, note that the proof of essential surjectivity of $K$ only requires $\phi$ to preserve the order.

	Moreover, by restricting the whole discussion to the commutative case, since $F\A$ is always commutative provided $\A$ is (\cite[Remark 5.3.4(iii)]{saito15monotone}), the same argument proves monadicity of the forgetful functors $\cSigmaC\to\cCalg $ and $\cSigmaCcpu\to \cCalgcpu $.
\end{remark}

We now seek to show that $\SigmaC$ is also cocomplete.
To this end, we show two preparatory lemmas, where we need to consider nonunital settings.
Recall from \cite[\S2.2]{saito15monotone} that a \newterm{nonunital \cstars{}} is a nonunital \cstar{}\footnote{Here, ``nonunital'' means ``not necessarily unital''; thus unital \cstar{}s are particular instances of nonunital \cstar{}s.} such that all norm-bounded monotone increasing sequences admit a supremum.
This is different from requiring that all monotone increasing sequences bounded by an element admit a supremum: this second option gives rise to \emph{pseudo monotone $\sigma$-complete (nonunital) \cstar{}s} and it allows one to consider the compact operators over a separable Hilbert space as an example, while the definition we adopt here does not (see \cite[Example 2.2.11]{saito15monotone} for details).

\begin{lemma}\label{lem:sigmaC_ideal}
	Let $\B$ be a \cstars{}, and let $I$ be a closed ideal in $\B$ and a nonunital \cstars{}.
	Then $I$ is a $\sigma$-ideal.
\end{lemma}
\begin{proof}
	Given a norm-bounded monotone increasing sequence $(a_n)\subseteq I$, we have $\sup_n^I a_n \ge \sup_n^{\B}a_n \ge a_n$ because $I\subseteq \B$ and both are \cstars{}s ($I$ nonunital). Since closed ideals are \emph{hereditary}~\cite[Theorem 3.2.1]{murphy90cstar}, that is, $0 \le b \le a$ with $a \in I$ and $b\in\B$ forces $b\in I$,
	we deduce that $\sup_n^{\B}a_n $ must also be contained in $I$, hence $\sup_n^{\B} a_n = \sup_n^I a_n$.
\end{proof}

\begin{lemma}\label{lem:reflexive_coeq}
	Consider \cstars{}s $\A$ and $\B$. Let us consider the diagram
	\[
		\xymatrix{
		\A
		\ar@<1.2ex>[r]^f
		\ar@<-1.2ex>[r]_g
		&
		\B
		\ar@{->}[l]|-{s}
		}
	\]
	where $f$ and $g$ are \starshom{}s, $s$ is any \cpu{} map, and $fs=gs=\id_{\B}$. The reflexive coequalizer $\C$ of the above diagram in $\Calgcpu$ is a \cstars{}, and moreover the induced $\B \to \C$ is a surjective \starshom{}.

	In other words, the forgetful functor $\SigmaC \to \Calgcpu$ preserves coequalizers that become reflexive after applying the functor.
	Consequently, the same holds for the forgetful functors $\SigmaC\to \Calg$ and $\SigmaC \to \SigmaCcpu$.
\end{lemma}
This statement is similar to \cite[Proposition 4.4]{furber2024monad}.
In particular, this result could be used to reprove monadicity of $U$ via Beck's theorem as stated in \cite[Theorem 3.14]{barr2005toposes}.

\begin{proof}
	The whole proof consists in proving that the $*$-subspace $I\coloneqq {\Im(f-g)}$ is a $\sigma$-ideal.
	Indeed, if this is the case, the coequalizer is given by the quotient map $\B \to \B/I$, and the required monotone properties follow immediately from \cite[Proposition 2.2.21]{saito15monotone}.

	First, we note that $I$ is a two-sided ideal.
	Indeed, given $a\in \A$ and $b\in \B$, then
	\[ b(f(a)-g(a)) = f(s(b)) f(a) - g(s(b)) g(a) = f(s(b)a)-g(s(b)a),\]
	and analogously on the other side.

	To show it is both closed and $\sigma$-closed, we consider the (nonunital) \starshom{}
	\[
		f'\colon
		\xymatrix{\ker(g)\, \ar@{^{(}->}[r] & \A \ar[r]^f& \B}
	\]
	and claim that its image is precisely $I$.
	Indeed, given $k \in \ker(g)$, $f'(k)=f(k)-g(k)\in I$.
	Conversely, if $b\in I$, there is $a\in \A$ such that $b=f(a)-g(a)$.
	Now, by selecting $a'=a-sg(a)$, the equations $g(a')=0$ and $f(a')=b$ are directly verified.
	Finally, $I=\Im(f')$.
	Consequently, $I$ is closed by \cite[Theorem 3.1.6]{murphy90cstar}.
	Moreover, using the first isomorphism theorem for \cstar{}s, $I = {\Im(f')}\cong \ker(g)/\ker(f')$, i.e.~it is a quotient of (nonunital) \cstars{}s, thus a nonunital \cstars{} by \cite[Proposition 2.2.21]{saito15monotone}.
	We have shown that $I\subseteq \B$ satisfies the hypotheses of \cref{lem:sigmaC_ideal}, so it is $\sigma$-closed.
\end{proof}

\begin{proposition}\label{prop:cocomplete}
	$T$ preserves reflexive coequalizers.
	In particular, $\SigmaC$ is cocomplete.
\end{proposition}
\begin{proof}
	Consider $f,g\colon \A\to\B$ and $s\colon \B \to \A$ with $gs=fs=\id_{\B}$ in $\Calg$.
	We call $\C$ and $\D$ the reflexive coequalizers of the pairs $(f,g)$ and $(Tf,Tg)$, respectively.
	By \cref{lem:reflexive_coeq}, $\D$ is a \cstars{} and moreover the cocone map $T\B \to \D$ is a \starshom{}.
	Precomposing with $\eta_{\B}\colon \B \to T\B$, we obtain a map $\B \to \D$ that coequalizes $f$ and $g$.
	This entails the existence of a unique map $\C \to \D$, which lifts to a map $T\C \to \D$ because $\D$ is a \cstars{}.
	Since the necessary commutativities are directly inferred, $T\C$ must also be the coequalizer.

	The second sentence follows from Linton's theorem (\cite[Corollaries 2 and 3]{linton1969coequalizers}), since $T$ is a monad on a cocomplete category (\cite[Section 2]{khoshkam2002categorical}).
\end{proof}
\begin{remark}\label{rem:two_routes_cocomplete}
		Cocompleteness of $\SigmaC$ admits a second, less elementary proof: by \cref{thm:U_accessible} the category is locally $\aleph_1$-presentable, hence cocomplete~\cite[Definition 1.17]{adamek1994locallypresentable}.
		We have given the argument above because it is independent of the accessibility machinery of \cref{sec:accessible} and because its first assertion, that $T$ preserves reflexive coequalizers, is needed later (it is the hypothesis of \cref{thm:SMC2} used in the proof of \cref{thm:tensor_shom}) and does not follow from presentability.
		No circularity arises: the proof of \cref{thm:U_accessible} rests on \cref{thm:U_monadic} and \cref{lem:technical_access}, not on \cref{prop:cocomplete}.
\end{remark}

\begin{remark}[Monadicity of (commutative) \wstar{}s]\label{rem:monadicity_W}
	Consider the category $\Walg$ of \wstar{}s and normal (unital) $*$-homomorphisms.
	The forgetful functor $\Walg \to \Calg$ has been proven to be monadic in \cite[Theorem 4.6]{furber2024monad}.
	Moreover, \cite[Proposition 4.4]{furber2024monad} also ensures that this functor preserves reflexive coequalizers.
	In particular, $\Walg$ is cocomplete by the same reasoning of the proof of \cref{prop:cocomplete}.

	Combining this result with \cref{thm:U_monadic} and \cite[Corollary 4.5.7]{borceux94handbook2}, we immediately obtain that $\Walg\to \SigmaC$ is monadic as well.

	A similar result holds by restricting to commutative algebras.
\end{remark}

We conclude with a simple remark that will be used in the proof of \cref{thm:tensor_scpu}.
\begin{remark}\label{rem:T_surjective}
	$T$ preserves surjective maps.
	Indeed, given a surjective map $\phi \colon \A \to \B$ in $\Calg$, the map $T\phi \colon T\A \to T\B$ sends $\eta_{\A}(a) $ to $\eta_{\B}\phi(a)$. 
	Since $\phi$ is surjective, $\Im (T\phi)\supseteq \eta_{\B}(\B)$. 
	As $T\phi$ is $\sigma$-normal, its image must be $\sigma$-closed, so it must coincide with the whole $T\B$ because $T\B$ is the $\sigma$-closure of $\eta_{\B}(\B)$.
\end{remark}

%% file: sections/04_giry_comonads.tex
\section{Noncommutative and abstract Giry (co)monads}\label{sec:giry}
As a consequence of monadicity, we can infer crucial results regarding the coKleisli nature of $\SigmaCcpu$ and $\cSigmaCcpu$ by virtue of the adjoint lifting theorem~\cite[Theorem 4.5.6]{borceux94handbook2}.

We start with a preliminary result, which previously appeared in the context of pre-\cstar{}s~\cite[Proposition 4.18]{fritz2023involutive}.
For the sake of transparency, we offer a complete argument.

Recall that, given a \cstar{} $\A$, the \newterm{state space} $S(\A)$ is the set of states $\phi\colon \A \to \bbC$ equipped with the weakest topology making the evaluations
\[
	\ev_a \colon \xymatrix@R=0.1em{
	S(\A) \ar[r]
	& \bbC \\
	\phi  \ar@{|->}[r]
	& \phi(a)
	}
\]
continuous for all $a\in \A$.
It is well-known that $S(\A)$ is a compact Hausdorff space: for instance, this follows from the discussion of \cite[3.2.1]{pedersen2018automgroups}, which shows that $S(\A)$ is a compact subspace of $A^*$ equipped with the weak-$*$ topology.

\begin{proposition}\label{prop:state_space_leftadjoint}
	The forgetful functor $R_c\colon \cCalg \to \Calgcpu$ has a left adjoint given by\footnote{The following proof strategy was suggested to us by Tobias Fritz.}
	\[
		L_c\colon \xymatrix@R=.1em{
		\Calgcpu \ar[r]
		& \cCalg \\
		\A  \ar@{|->}[r]
		& C(S(\A))
		}
	\]
\end{proposition}
\begin{proof}
	Let $\CHaus$ be the category of compact Hausdorff spaces with continuous maps.
	For brevity, maps $\phi\colon \A \times X\to \bbC$, where $\A$ is a \cstar{} and $X$ is a compact Hausdorff space, are called \emph{state-continuous maps} if $\phi(-,x)\colon \A \to \bbC$ are states for all $x \in X$ and $\phi(a,-)\colon X\to \bbC$ are continuous for all $a\in \A$  (cf.\ the definition of Markov kernels, recalled on p.~\pageref{def:mkernel}); the set of these maps will be denoted by $\mathbf{SCM}(\A\times X,\bbC)$.

	The statement is obtained by the following chain of isomorphisms (whose naturality in $\A$ and $X$ is a direct check):
	\begin{align*}
		\cCalg (C(S(\A)),C(X)) & \cong \CHaus(X,S(\A))                                        \\
		                       & \cong \mathbf{SCM}(\A\times X,\bbC) \tag{$\spadesuit$} \\
		                       & \cong \Calgcpu(\A,C(X)) \tag{$\heartsuit$}
	\end{align*}
	While the first isomorphism is an immediate consequence of Gelfand duality, the remaining two require further explanation.
	We do so by exhibiting mutually inverse maps.
	\begin{itemize}
		\item For $(\spadesuit)$, we define
		      \[\begin{gathered}
			      \xymatrix@R=.1pt{
			      \CHaus(X,S(\A)) \ar[r]
			      & \mathbf{SCM}(\A\times X,\bbC) \\
			      \phi  \ar@{|->}[r]
			      & \bar{\phi}(a,x)\coloneqq \phi(x)(a)
			      } \\
			      \xymatrix@R=.1pt{
			      \mathbf{SCM}(\A\times X,\bbC)  \ar[r]
			      & \CHaus(X,S(\A)) \\
			      \psi  \ar@{|->}[r]
			      & \bar{\psi}(x)\coloneqq  \psi(-,x)
			      }
								\end{gathered}\]
		      These maps are well-defined: $\bar{\phi}$ is state-continuous because $\phi(x)$ is a state by definition of $S(\A)$ and $\phi(-)(a)=\ev_a\circ\, \phi$ is continuous because it is a composition of continuous maps. For the second one, by definition of $S(\A)$, continuity of $\bar{\psi}$ can be checked by postcomposition with $\ev_a$, and $\ev_a \circ \bar{\psi} = \psi(a,-)$ is continuous by assumption.
		      That these are inverses of each other is a direct check.

		\item For $(\heartsuit)$, we consider
		      \[\begin{gathered}
			      \xymatrix@R=.1pt{
										\mathbf{SCM}(\A\times X,\bbC)  \ar[r]
										& \Calgcpu(\A,C(X)) \\
										\psi  \ar@{|->}[r]
										& \widetilde{\psi}(a)\coloneqq  \psi(a,-)
										} \\
										\xymatrix@R=.1pt{
											\Calgcpu(\A,C(X)) \ar[r]
											& \mathbf{SCM}(\A\times X,\bbC) \\
											\phi  \ar@{|->}[r]
											& \widetilde{\phi}(a,x)\coloneqq \phi(a)(x)
											}
										\end{gathered}
		      \]
		      Concerning the well-definedness of $\widetilde{\psi}$, we have that $\psi(a,-)$ is continuous by definition, so $\widetilde{\psi}(a)$ indeed belongs to $C(X)$.
		      To show that $\widetilde{\psi}$ is a \cpu{} map, take any positive element $a\in \A$: since $\psi(-,x)$ is a state, $\psi(a,x)$ is nonnegative.
		      As this holds for all $x$, $\psi(a,-)$ is a positive element of $C(X)$.
		      A similar reasoning ensures that $\widetilde{\psi}(1)=1$, and therefore $\widetilde{\psi}$ is a \cpu{} map because its codomain is commutative (\cite[Theorem 3.9]{paulsen2002}).
		      Finally, $\widetilde{\phi}$ is a state-continuous map because $\widetilde{\phi}(a,-)$ is continuous by definition, and $\widetilde{\phi}(-,x)$ is a state for $\A$ because it is a composition of $\phi$ with the state $C(X)\owns f\mapsto f(x)\in \bbC$.
		      As above, these maps are inverses of one another by direct check.
		      \qedhere
	\end{itemize}
\end{proof}

\begin{theorem}[Noncommutative and abstract Giry (co)monads]\label{thm:comonad}
The following assertions hold.
	\begin{enumerate}
		\item\label{it:SigmaCcpu_coKleisli} The forgetful functor $R'\colon \SigmaC\to \SigmaCcpu$ has a left adjoint $L'$, and the induced comonad $L'R'$ on $\SigmaC$ has $\SigmaCcpu$ as its coKleisli category. This also holds when restricting to Pedersen--Baire envelopes.
		\item\label{it:cSigmaCcpu_coKleisli} Analogously, the forgetful $R'_{cc}\colon \cSigmaC\to \cSigmaCcpu$ has a left adjoint $L'_{cc}$, and the induced comonad $L'_{cc}R'_{cc}$ on $\cSigmaC$ has $\cSigmaCcpu$ as its coKleisli. This also holds when restricting to commutative Pedersen--Baire envelopes.
		\item\label{it:state_space} The forgetful functor $R'_c\colon \cSigmaC \to \SigmaCcpu$ has a left adjoint. This also restricts to Pedersen--Baire envelopes.
	\end{enumerate}
\end{theorem}
\begin{proof}
	We focus on Item~\ref{it:SigmaCcpu_coKleisli}.
	By \cref{thm:U_monadic,rem:monadic_cpu}, the forgetful functors $U\colon \SigmaC \to \Calg$ and $\bar{U}\colon \SigmaCcpu \to \Calgcpu$ are monadic.
	Moreover,
	\begin{itemize}
		\item The forgetful functor $R\colon \Calg \to \Calgcpu$ has a left adjoint $L\colon \Calgcpu \to \Calg$ by \cite[Proposition 8]{kirchberg1998Cstarfromoperator}.
		\item The forgetful functor $R'\colon \SigmaC \to \SigmaCcpu$ satisfies $\bar{U}\circ R' = R\circ U$ by immediate check (this composition is the obvious forgetful functor).
	\end{itemize}
	Putting it all together, we obtain the following diagram
	\[
		\vxy{
		\SigmaC\ar[r]^-{R'} \ar[d]^U & \SigmaCcpu \ar[d]^{\bar{U}}\\
		\Calg \ar@<.4em>[r]^-R\ar@{<-}@<-.4em>[r]_-L\ar@{}[r]|-\top & \Calgcpu}
	\]
	whose vertical arrows are monadic.
	By the adjoint lifting theorem~\cite[Theorem 4.5.6]{borceux94handbook2}, since $\SigmaC$ is cocomplete by \cref{prop:cocomplete}, we conclude that $R'$ has a left adjoint $L'$.
	The fact that the coKleisli category of the comonad $L' \circ R'$ is isomorphic to $\SigmaCcpu$ follows from the fact that $R'$ is bijective on objects: for all \cstars{}s $\A$ and $\B$,
	\begin{align*}
		\coKl(L'R') (\A,\B) & \cong \SigmaC(L'R'(\A) , \B)\cong \SigmaCcpu(R'(\A),R'(\B)) \\
		                    & \cong \SigmaCcpu(\A,\B).
	\end{align*}
	The restriction to Pedersen--Baire envelopes is straightforward, as the diagram of left adjoint functors commutes because the diagram of right adjoints does: for all $\A$ and $\B$ in $\SigmaC$,
	\begin{align*}
		\SigmaC(L'\bar{F} \A,\B) & \cong  \SigmaCcpu (\bar{F}\A, R'\B) \cong \Calgcpu (\A,\bar{U}R' \B)      \\
		                         & \cong \Calgcpu (\A, RU\B) \cong \Calg(L\A,U\B) \cong \SigmaC(FL \A , \B),
	\end{align*}
	and so by the Yoneda lemma $L'\bar{F} \A\cong FL \A$.

	The proof of Items~\ref{it:cSigmaCcpu_coKleisli} and \ref{it:state_space} are analogous.
	They are obtained by the commutative diagrams
	\[
		\vxy{
		\cSigmaC\ar@<.4em>[r]^-{R'_{cc}}\ar@{<-}@<-.4em>[r]_-{L'_{cc}}\ar@{}[r]|-\top \ar[d]^U & \cSigmaCcpu \ar[d]^{\bar{U}}\\
		\cCalg \ar@<.4em>[r]^-{R_{cc}}\ar@{<-}@<-.4em>[r]_-{L_{cc}}\ar@{}[r]|-\top & \cCalgcpu} \qquad \vxy{
		\cSigmaC\ar@<.4em>[r]^-{R'_c}\ar@{<-}@<-.4em>[r]_-{L'_c}\ar@{}[r]|-\top \ar[d]^U & \SigmaCcpu \ar[d]^{\bar{U}}\\
		\cCalg \ar@<.4em>[r]^-{R_c}\ar@{<-}@<-.4em>[r]_-{L_c}\ar@{}[r]|-\top & \Calgcpu}
	\]
	where $U$ and $\bar{U}$ in both diagrams are monadic as observed in \cref{rem:monadic_cpu}, the left adjoints of $R_{cc}$ and $R_c$ are given by the state space (see \cref{prop:state_space_leftadjoint}), and therefore $L'_{cc}$ and $L'_c$ are again obtained by the adjoint lifting theorem.\footnote{Cocompleteness of $\cSigmaC$ can be obtained by adapting our \cref{prop:cocomplete} above. Alternatively, this holds because $\cSigmaC$ is equivalent to the category of \bsigma{}s and $\sigma$-homomorphisms, which is known to be complete and cocomplete (\cite[Remark 2.1.12]{fritz2025categories}). }
	It is worth stressing that while $R'_{cc}$ is bijective on objects, $R'_c$ is not: for this reason, $\cSigmaCcpu$ is a coKleisli category over $\cSigmaC$, whereas the same does not hold for $\SigmaCcpu$.
	The restriction to (commutative) Pedersen--Baire envelopes holds as before.
\end{proof}
	The name given to the previous result is not merely suggestive.
	We show in \cref{thm:giry_restriction} that, under the equivalence between Baire measurable spaces and commutative Pedersen--Baire envelopes, the comonad of \cref{thm:comonad} \emph{is} the Giry monad.
	The remainder of this section makes this precise, and may be skipped by the reader willing to take the name on trust.

The Giry monad, first appeared in \cite{giry82probability}, is a monad $\mathcal{G}$ acting on $\Meas$, the category of measurable spaces with measurable functions.
Given a measurable space $X$, $\mathcal{G}(X)$ is given by the set of probability measures on $X$, equipped with the smallest $\sigma$-algebra making
\[
	\ev_U \colon \xymatrix@R=0pt{
	\mathcal{G}(X) \ar[r]
	& \bbC \\
	p  \ar@{|->}[r]
	& p(U)
	}
\]
measurable for all measurable subsets $U\subseteq X$ (cf.\ the state space of a \cstar{}).

The Kleisli category associated to $\mathcal{G}$ is $\Stoch$, the category of measurable spaces with Markov kernels (\cite[Section I.4]{giry82probability}).
As discussed in \cref{sec:preliminaries}, \scpu{} maps between commutative \cstars{}s may be viewed as ``abstract Markov kernels'', and similarly, one may view \scpu{} maps between \cstars{}s as ``noncommutative Markov kernels'', thus establishing a connection with the previous statement.
For the sake of convenience, we recall that a \newterm{Markov kernel} from $(X,\Sigma_X)$ to $(Y,\Sigma_Y)$ is a map $\phi\colon \Sigma_Y \times X \to [0,1]$ such that $\phi(U,-)\colon X\to [0,1]$ is measurable for all $U\in \Sigma_Y$ and $\phi(-,x)\colon \Sigma_Y\to [0,1]$ is a probability measure for all $x\in X$.\label{def:mkernel}

As argued in \cite{jamneshan23foundational}, a natural framework for measure theory that avoids second-countability constraints, unlike standard Borel spaces, is provided by \newterm{Baire measurable spaces}. These are measurable spaces $(X,\Sigma_X)$ where $\Sigma_X$ is the Baire $\sigma$-algebra of a compact Hausdorff topology on $X$.
We now seek to show that the Giry monad, when restricted to Baire measurable spaces, in fact coincides with the comonad structure obtained in \cref{thm:comonad}. 

\begin{theorem}\label{thm:giry_restriction}
	Let us consider:
	\begin{itemize}
		\item The category $\BaireMeas$ of Baire measurable spaces with measurable functions.
		\item The category $\cPBenv$ of commutative Pedersen--Baire envelopes with \starshom{}s.
	\end{itemize}
	Then the Giry monad restricts to $\BaireMeas$, i.e.\ $\mathcal{G}(X)$ is Baire measurable whenever $X$ is, and the following diagram commutes up to a natural isomorphism
	\[
		\vxy{
		\BaireMeas^{\op}\ar[r]^-{\mathcal{G}}\ar[d]^{\Linf} & \BaireMeas^{\op} \ar[d]^{\Linf}\\
		\cPBenv \ar[r]^-{\mathcal{S}'} & \cPBenv}
	\]
	where $\mathcal{S}'$ denotes the (restricted) comonad obtained by the forgetful functor $R'_{cc}\colon \cSigmaC \to \cSigmaCcpu$ and its adjoint.

\end{theorem}
The vertical arrows in the statement are equivalences \cite[Corollary 4.2.9]{fritz2025categories}.
\begin{proof}
	The fact that $\mathcal{G}(X)$ is Baire measurable whenever $X$ is follows from \cite[Corollary A.4]{shirazi2026commutativity}.
	To show the desired commutativity, we argue with the respective characterizations of the Kleisli and coKleisli categories.
	\begin{itemize}
		\item When restricted to $\BaireMeas$, the Kleisli category of the Giry monad is $\BaireStoch$, the category of Baire measurable spaces with Markov kernels. This follows from the fact that the Kleisli category of the Giry monad on $\Meas$ is $\Stoch$ (\cite[Section I.4]{giry82probability}).
		\item By \cref{thm:comonad}, the coKleisli category of $\mathcal{S}'$ is $\cpbenvcpu$, the category of commutative Pedersen--Baire envelopes with \scpu{} maps.
	\end{itemize}
	We claim that the diagram
	\begin{equation}\label{eq:Giry_extend_forget}
		\vxy{
			\BaireMeas^{\op}\ar[r]\ar[d]^{\Linf} & \BaireStoch^{\op} \ar[d]^{\Linf}\\
			\cPBenv \ar[r] & \cpbenvcpu}
	\end{equation}
	where the horizontal arrows are forgetful functors and the vertical arrows are equivalences \cite[Corollaries 4.2.9 and 4.3.8]{fritz2025categories}, is commutative.
	This is easily seen once all details have been made explicit.

	Let $\phi\colon X\to Y$ be an arbitrary measurable function in $\BaireMeas$. Via the forgetful $\BaireMeas\to \BaireStoch$, $\phi$ is sent to
	\[
		\bar{\phi} \colon \xymatrix@R=.1pt{
		\Sigma_Y\times X \ar[r]
		& [0,1] \\
		(U,x)
		\ar@{|->}[r]
		& \delta_{\phi(x)}(U),}
	\]
	where $\Sigma_Y$ is the $\sigma$-algebra of $Y$ and $\delta_{p}$ is the delta measure concentrated in $p$ (the unit of the Giry monad $X\to \mathcal{G}(X)$ is given by sending each point to its delta measure).
	The commutativity of \eqref{eq:Giry_extend_forget} is then equivalent to showing that $\Linf(\phi)= \Linf(\bar{\phi})$,
	which means that for every $f\colon Y\to \bbC$ and $x \in X$, $f(\phi(x))=\int_Y f(y) \bar{\phi}(\mathrm{d}y,x)$.

	But this holds by definition of $\bar{\phi}$:
	\[
		\int_Y f(y) \bar{\phi}(\mathrm{d}y,x)=\int_Y f(y) \delta_{\phi(x)}(\mathrm{d}y) =f(\phi(x)),
	\]
	so commutativity of \eqref{eq:Giry_extend_forget} is finally ensured.

	Since a square of left adjoints commutes whenever the corresponding square of right adjoints does, (as mentioned in the proof of \cref{thm:comonad}), and since $\Linf$ is an equivalence, we obtain the commutative diagram
	\[
		\vxy{
			\BaireStoch^{\op} \ar[r]\ar[d]^{\Linf} & \BaireMeas^{\op} \ar[d]^{\Linf}\\
			\cpbenvcpu \ar[r] & \cPBenv}
	\]
	where the horizontal arrows send $X$ to $\mathcal{G}(X)$ and $\A$ to $\mathcal{S}'(\A)$, respectively.
	We can now horizontally compose the two commutative diagrams to conclude the proof of the statement.
\end{proof}

\begin{remark}
	The categorically inclined reader may note that the proof above, making use of the uniqueness of left adjoints, directly implies that $\zeta\colon \Linf \mathcal{G}\Rightarrow \mathcal{S}'\Linf$ is compatible with the counit and comultiplication structures. Equivalently, the following diagrams of natural transformations between composite functors commute:
	\[\small
		\xymatrix{
		\Linf \circ\mathcal{G} \ar@{=>}[rr]^{\zeta} \ar@{=>}[dr]_{\Linf * \epsilon^{\mathcal{G}}} &&\mathcal{S}'\circ\Linf \ar@{=>}[dl]^{\epsilon^{\mathcal{S}'}* \Linf} \\
		& \Linf
		}% \qquad
		\xymatrix{
		\Linf \circ\mathcal{G} \ar@{=>}[rr]^{\zeta} \ar@{=>}[d]_{\Linf * \mu^{\mathcal{G}}}
		&& \mathcal{S}' \circ\Linf \ar@{=>}[d]^-{\mu^{\mathcal{S}'} * \Linf}
		& \\
		\Linf \circ\mathcal{G}\circ\mathcal{G} \ar@{=>}[r]_{\zeta *\mathcal{G}}
		& \mathcal{S}'\circ\Linf \circ\mathcal{G}\ar@{=>}[r]_{\mathcal{S}' *\zeta} & \mathcal{S}'\circ\mathcal{S}'\circ \Linf
		}
	\]
	where $\epsilon^{-} $ and $\mu^{-}$ denote the counit and comultiplication of the comonads $\mathcal{G}$ (on $\BaireMeas^{\op}$) and ${\mathcal{S}'}$.
\end{remark}

%% file: sections/05_local_presentability.tex
\section{Local \texorpdfstring{$\aleph_1$}{א{1}}-presentability of \texorpdfstring{\cstars{}s}{σC*-algebras}}\label{sec:accessible}

In this section, we ensure local $\aleph_1$-presentability of $\SigmaC$ and accessibility of the forgetful functor $U$.
The reader is advised to have some background on presentable and accessible categories (see, e.g., \cite[Section 1 and 2]{adamek1994locallypresentable} for some definitions that will be used below).
For the sake of readability, we recall that a poset $P$ is \newterm{$\aleph_1$-directed} if every countable subset has an upper bound.

We start with a technical result.

\begin{lemma}\label{lem:technical_access}
	Let us consider an $\aleph_1$-directed diagram $D\colon P \to \Calg$, by which we mean a functor whose domain $P$ is an $\aleph_1$-directed poset. Call $\A\coloneqq \colim D$, and $\phi_{\alpha',\alpha}\colon D(\alpha')\to D(\alpha)$ the transition morphisms, and $\phi_{\alpha}\colon D(\alpha)\to \A$ the maps of the initial cocone. Then,
	\begin{enumerate}
		\item\label{it:technical_access_norm}\emph{Isometric representation}. For every $a\in \A$, there exist $\alpha \in P$ and $a_{\alpha}\in D(\alpha)$ such that $a=\phi_{\alpha}(a_{\alpha})$.
		      Moreover, $a_{\alpha}$ can be chosen so that $\norm{a_{\alpha}}_{D(\alpha)} = \norm{a}_{\A}$.
		\item\label{it:technical_access_positivity}\emph{Positive representation}. For every positive $a\in \A$, there exists $\alpha\in P$ and a positive $a_{\alpha}\in D(\alpha)$ such that $a=\phi_{\alpha}(a_{\alpha})$.
	\end{enumerate}
\end{lemma}
\begin{proof}\item
	\begin{enumerate}
		\item[\ref{it:technical_access_norm}]
		      Consider the quotient $B=\bigcup_{\alpha \in P} D(\alpha)/_\sim$, i.e.\@\xspace the colimit of the diagram $D$ in the category of $*$-algebras.
		      As discussed in \cite[Corollaries 7.22 and 7.23]{bunke2021lectures}, $\A$ corresponds to the norm-completion of $B$. %(a fact we already used in \cref{prop:tensor_colims}).     

		      We claim that for $a\in B$, we can find $\alpha$ and $a_{\alpha}\in D(\alpha)$ such that $\phi_{\alpha}(a_{\alpha})=a$ and $\norm{a_{\alpha}}_{D(\alpha)}=\norm{a}_B$.
		      Recall that the norm is defined as
		      \[\norm{a}_B \coloneqq \inf \norm{a_{\alpha}}_{D(\alpha)},\]
		      where the infimum is calculated over all $\alpha \in P$ and $a_{\alpha}\in D(\alpha)$ such that $\phi_{\alpha}(a_{\alpha})=a$.
		      In particular, for each $n$, we can find $\alpha_n$ and $a_n\in D(\alpha_n)$ such that $\phi_{\alpha_n}(a_n)=a$ and $\norm{a_n}_{D(\alpha_n)} \le \norm{a}_B+ \frac{1}{n}$.
		      By definition of $B$, $a_n \sim a_{m}$ for all $n,m$.
		      Since $P$ is directed, this means that there exists a certain $\alpha>\alpha_n,\alpha_m$ such that $\phi_{\alpha_n,\alpha} (a_n) = \phi_{\alpha_m,\alpha}(a_m)$.
		      By letting $m$ vary, since every countable subset of $P$ has an upper bound, such an $\alpha$ can be chosen to identify all the elements $a_n$, i.e.\  $\phi_{\alpha_n,\alpha} (a_n) = \phi_{\alpha_m,\alpha}(a_m)$ for all $m,n$ simultaneously.
		      We now define $a_{\alpha}\coloneqq \phi_{\alpha_n,\alpha} (a_n)\in D(\alpha)$, which does not depend on $n$ by the argument above.
		      Since $\phi_{\alpha_n,\alpha}$ is a $*$-homomorphism, it is contractive, and therefore
		      \[ \norm{a_{\alpha}}_{D(\alpha)}\le \norm{a_n}_{D(\alpha_n)}\le \norm{a}_B+ \frac{1}{n} \]
		      from which $\norm{a_{\alpha}}_{D(\alpha)}= \norm{a}_B$ by arbitrariness of $n$.
		      Finally, the claim is proved.

		      We now seek to show that $B=\A$.
		      Since $\A$ corresponds to the norm-completion of $B$ (as mentioned at the beginning of the proof), it suffices to take any Cauchy sequence $(a_n)$ of $B$ and show that it admits a limit in $B$.

		      Passing to a subsequence if necessary, we may assume that $\norm{a_n-a_{n-1}}\le \frac{1}{2^{n}}$ for all $n\ge 2$.
		      For each difference $d_n\coloneqq a_n-a_{n-1}$, with $n\ge 1$ (we set $a_0\coloneqq 0$), one can find $\alpha_n$ and $d'_n$ such that $\phi_{\alpha_n}(d'_n)=d_n$ and $\norm{d'_n}=\norm{d_n}$, hence $\norm{d'_n}\le \frac{1}{2^{n}}$.
		      By taking the supremum $\alpha$ of all $\alpha_n$, we can assume that $(d'_n)\subseteq D(\alpha)$.
		      In other words, we identify $d'_n$ with its image in $D(\alpha)$.
		      We note that the norms are not increased by this assumption, since all transition morphisms $\phi_{\alpha_n,\alpha}$ are contractive.
		      Then, the sequence $s_n \coloneqq \sum_{i=1}^n d'_i$ is Cauchy in $D(\alpha)$ by direct check: for $m>n\ge 1$,
		      \[
			      \norm{s_m - s_n} =\norm{\sum_{i=n+1}^m d_i} \le \sum_{i=n+1}^m \frac{1}{2^i} = \frac{1}{2^n} \left( \sum_{i=1}^{m-n} \frac{1}{2^i}\right) \le \frac{1}{2^n} \xrightarrow{n,m\to \infty} 0.
		      \]
		      Moreover, its limit $a_{\alpha}$ in $D(\alpha)$ yields a limit for $(a_n)$ because $\phi_{\alpha}$ is bounded:
		      \[
			      \phi_{\alpha}(a_{\alpha}) =\lim_n \phi_{\alpha}(s_n) = \lim_n \left(\sum_{i=1}^n d_i \right)= \lim_n a_n.
		      \]
		      Finally, $B=\A$, so the statement follows by the claim we proved above.
		\item[\ref{it:technical_access_positivity}] Let $a\in \A$ be a positive element.
		      By Item~\ref{it:technical_access_norm}, there exists some $\alpha$ and $a'_{\alpha}$ satisfying $\phi_{\alpha}(a'_{\alpha})=a$.
		      We now take $a_{\alpha}\coloneqq \max(a'_{\alpha},0)\ge 0$, where this notation is a standard application of the continuous functional calculus (see, e.g., \cite[Proposition 4.6]{takesaki1979operatoralgebras1}).
		      By functoriality of functional calculus~\cite[Corollary 8.4]{takesaki1979operatoralgebras1},
		      \[
			      \phi_{\alpha}(a_{\alpha})= \phi_{\alpha}(\max(a'_{\alpha},0))=\max(\phi_{\alpha}(a'_{\alpha}),0)=\max(a,0) = a
		      \]
		      where the last equality follows from $a\ge 0$.\qedhere
	\end{enumerate}
\end{proof}

\begin{remark}\label{rem:Calg_accessible}
	We recall that $\Calg$ is a locally $\aleph_1$-presentable category (combine \cite[Theorem 2.4]{pelletier1989equationalCstar} and \cite[Theorem 3.28]{adamek1994locallypresentable}). %\url{https://arxiv.org/pdf/2203.12087}. 
    The result of Pelletier and Rosick\'y quoted here is the closest antecedent of the present programme: it exhibits \cstar{}s as the models of an equational theory. 
	We establish a weaker analogue for \cstars{}s, namely that they are the models of an $\aleph_1$-ary essentially algebraic theory (\cref{cor:essentially_algebraic}). 
	Whether they are also the models of an equational theory is left open.
\end{remark}

\begin{theorem}\label{thm:U_accessible}
	The functor $U\colon \SigmaC \to \Calg$ is $\aleph_1$-accessible and $\SigmaC$ is locally $\aleph_1$-presentable.
\end{theorem}
\begin{proof}
	Let us consider an $\aleph_1$-directed diagram $D\colon P \to \SigmaC$, where $P$ is an $\aleph_1$-directed poset.
	We will show that $U(\colim D)\cong \colim UD$. 
	By \cite[Theorem 1.5 and Remark 1.21]{adamek1994locallypresentable} this will be enough to ensure a similar result holds for all colimits over $\aleph_1$-\emph{filtered} categories, not necessarily posets. 

	Concretely, since $U$ simply forgets $\sigma$-normality, it suffices to ensure that $\A \coloneqq \colim UD$ is a \cstars{} and the cocone morphisms are $\sigma$-normal.
	Once this is established, since left adjoints always preserve colimits, it is clear that $UF=T\colon \Calg \to \Calg$ is an $\aleph_1$-accessible monad over a locally $\aleph_1$-presentable category (see \cref{rem:Calg_accessible}), so $\SigmaC$ is also locally $\aleph_1$-presentable by combining \cref{thm:U_monadic} and \cite[Theorem and Remark 2.78]{adamek1994locallypresentable}.
	Finally, $U$ is indeed $\aleph_1$-accessible as claimed.

	By \cref{lem:technical_access}\ref{it:technical_access_norm}, all elements $a\in \A\coloneqq\colim UD$ are given by some $a_{\alpha}\in D(\alpha)$, where $\alpha$ is arbitrary.
	From now on, let $(a_n)$ be a norm-bounded monotone increasing sequence in $\A$. %We aim to show that $\sup_n a_n$ exists in $\A$.
	Below, we use $\phi_{\alpha',\alpha}$ and $\phi_{\alpha}$ to denote the transition morphisms and the cocone morphisms, respectively, as in the statement of \cref{lem:technical_access}.

	We first show that there exists $\gamma \in P$ and a norm-bounded monotone increasing sequence $(b_n)\subseteq D(\gamma)$ satisfying $\phi_{\gamma}(b_n)=a_n$.
	Let us take $\alpha_1\in P$ such that $a_1$ has a representative $d_1\in D(\alpha_1)$.
	Next, for $n>1$, note that $a_n-a_{n-1}\ge 0$ by assumption.
	Then there exists $\alpha_n\in P$ and a positive representative $d_n \in D(\alpha_n)$ of ${a_n - a_{n-1}}$ by \cref{lem:technical_access}\ref{it:technical_access_positivity}.

	By taking the supremum $\alpha_{\bbN}$ of $(\alpha_n)_n$, $a_n$ has a representative $b_n$ in $D(\alpha_{\bbN})$ given by $b_n\coloneqq \sum_{j=1}^n d_n$.
	By definition, $b_n \le b_{n+1}$.
	Moreover, by \cref{lem:technical_access}\ref{it:technical_access_norm}, for each $a_n$ one can find $\beta_n$ and $c_n \in D(\beta_n)$ such that $\phi_{\beta_n}(c_n) =a_n$ and $\norm{c_n}=a_n$.
	Since $c_n \sim b_n$, we can find a certain $\gamma_n > \alpha_{\bbN} , \beta_n $ such that $\phi_{\alpha_{\bbN},\gamma_n}(b_n) = \phi_{\beta_n, \gamma_n}(c_n)$; in particular, $\norm{\phi_{\alpha_{\bbN},\gamma_n}(b_n)}=\norm{a_n}$ because the transition morphisms are contractive.
	By now taking the upper bound $\gamma$ of all $\gamma_n$, and by an abuse of notation calling $b_n$ the elements $\phi_{\alpha_n ,\gamma}(b_n)\in D(\gamma)$, we can finally ensure that $(b_n)\subseteq D(\gamma)$ is a norm-bounded monotone increasing sequence.

	Next, given $x \ge a_n$ and a norm-bounded monotone increasing sequence $(b_n)\subseteq D(\gamma)$ satisfying $\phi_{\gamma}(b_n)=a_n$,
	we claim that there exists $\alpha>\gamma$ and a representative $y\in D(\alpha)$ of $x$ such that $\phi_{\gamma,\alpha}(\sup_n b_n)\le y$.
	This immediately ensures that $\phi_{\gamma}(\sup_n b_n)$ is indeed the supremum of $a_n$ in $\A$, so $\A$ is a \cstars{}, and moreover the arbitrariness of $(b_n)$ implies that the cocone morphisms are $\sigma$-normal.
	Therefore, this claim suffices to conclude the proof of the statement.

	We are now left to prove the claim of the previous paragraph. By \cref{lem:technical_access}\ref{it:technical_access_positivity}, $x-a_n\ge 0$ has a positive representation in some $D(\delta_n)$. By taking the upper bound $\delta$ of $(\delta_n)_n$ and $\gamma$,
	we have a representative $y_n\ge 0$ of $x-a_n$ in $D(\delta)$ for all $n$.
	Since $\delta\ge \gamma$, $y_n + \phi_{\gamma,\delta}(b_n) \sim y_{m}+ \phi_{\gamma,\delta}(b_m)$ because they are both sent to $x$ via the cocone morphism.
	By letting $m$ vary, we can find some $\alpha\ge \delta$ such that $y\coloneqq \phi_{\delta,\alpha}(y_n + \phi_{\gamma,\delta}(b_n)) = \phi_{\delta, \alpha}(y_{m}+ \phi_{\gamma,\delta}(b_m))$ for all $n,m$.
	Then $y\ge \phi_{\gamma,\alpha}(b_n)$ for all $n$ because $\phi_{\delta,\alpha}(y_n)\ge 0$, and thus $y\ge \sup_n \phi_{\gamma,\alpha}(b_n) = \phi_{\gamma,\alpha}(\sup_n b_n)$, because transition morphisms are $\sigma$-normal.
\end{proof}

\begin{corollary}\label{cor:essentially_algebraic}
	$\SigmaC$ is the category of models of an $\aleph_1$-ary essentially algebraic theory.
\end{corollary}
\begin{proof}
		A category is locally $\lambda$-presentable if and only if it is the category of models of a $\lambda$-ary essentially algebraic theory~\cite[Theorem 3.36]{adamek1994locallypresentable}; it therefore suffices to apply \cref{thm:U_accessible}.
\end{proof}

\begin{remark}\label{rem:presentability_consequences}
\cref{thm:U_accessible} carries further consequences. 
For instance, $\SigmaC$ is wellpowered and co-wellpowered. 
Moreover, the family of Pedersen--Baire envelopes $F\A$ of $\aleph_1$-presentable \cstar{}s $\A$ forms a strong generator of $\SigmaC$, and every \cstars{} is an $\aleph_1$-directed colimit of $\aleph_1$-presentable \cstars{}s.
The latter assertion is the precise form of the informal statement, made in the introduction, that each \cstars{} is built out of a small family of $\aleph_1$-compact ones.
\end{remark}

%% file: sections/06_tensor_product.tex
\section{A tensor product of \texorpdfstring{\cstars{}s}{σC*-algebras}}\label{sec:tensor}

Our next aim is to ensure that $\SigmaC$ and $\SigmaCcpu$ have a tensor product.
From now on,
\[
	\otimes\colon \Calg \times \Calg \to \Calg
\]
denotes the \newterm{maximal tensor product}, satisfying the following universal property: for every $\phi \colon \A \to \C$ and $\psi \colon \B \to \C$ with commuting images, there exists a unique $(\phi,\psi)\colon \A \otimes \B \to \C$ such that $(\phi,\psi)(a\otimes b)= \phi(a)\psi(b)$.
The resulting triple $(\Calg,\otimes,\bbC)$ is a symmetric monoidal category.
Throughout, we shall work exclusively with this monoidal structure.\footnote{In particular, the minimal tensor product will not be discussed.}

This tensor is known to extend to cpu maps with the same universal property, so in particular also $(\Calgcpu,\otimes ,\bbC)$ is a symmetric monoidal category.
For more details, the reader may refer to \cite{brownozawa}, particularly Section 3.3 and Exercise 3.5.1.

To equip $\SigmaC$ with a tensor product, we rely on its monadicity over $\Calg$ (\cref{thm:U_monadic}). 
The construction follows the well-known description of tensor products of algebras over monoidal monads, first investigated by Kock in the early 1970s~\cite{Kock1970}.
We employ results from \cite{seal2013tensors}, where this machinery was generalized to monoidal categories that are not necessarily closed, as is the case for $(\Calg,\otimes,\bbC)$.

Since monoidality of the monad is necessary, we focus on this first.
To prove it, we only rely on well-established results, so the reader is not required to know the precise definition.
\begin{proposition}\label{prop:monoidalmonad}
	The monad $T=UF$ obtained from the adjunction of \autoref{thm:adj_sigma} is a symmetric monoidal monad in the sense of \cite{Kock1970}.
\end{proposition}

\begin{proof}
	The Kleisli category of $T$ is given by Pedersen--Baire envelopes with $\sigma$-homomorphisms by direct check.
	This was proven to be symmetric monoidal in \cite[Proposition 5.3.2]{fritz2025categories}, and so $T$ is a symmetric monoidal monad by \cite[Proposition 1.2.2]{seal2013tensors}.
\end{proof}
Monoidality of $T$ is witnessed by the maps
\[ \theta_{\A,\B}\colon T\A \otimes T\B \to T(\A\otimes \B),\]
uniquely determined by the universal property of the maximal tensor product applied to the $\sigma$-homomorphisms $T(\A)\to T(\A\otimes B)$ and $T(\B)\to T(\A\otimes \B)$, which are uniquely obtained by
\begin{equation*}%\label{eq:description_theta}
	\A\owns a\mapsto \eta_{A\otimes B}(a\otimes 1) \in T(\A\otimes B)\qquad \text{and}\qquad \B \owns b\mapsto \eta_{A\otimes B}(1\otimes b) \in T(\A\otimes B).
\end{equation*}

% It is worth noting that, for order convergent sequences $(\eta_{\A} (a_n))_n\subseteq T\A$ and $(\eta_{\B} (b_m))_m\subseteq T\B$,\footnote{The use of order limits here is motivated by the fact that even when $(a_n)$ and $(b_m)$ are monotone increasing, $a_n \otimes b_m$ may not be. For example, if $b_m$ is negative, then $(a_n \otimes b_m)_n$ is monotone decreasing.} 
% \begin{equation}\label{eq:theta_description}
%   \theta_{\A,\B}(
%   \olim_n^{T\A}\eta_{\A}(a_n)\otimes
%   \olim_m^{T\B}\eta_{\B}(b_m)) =  \olim_n^{T(\A\otimes\B)}
%       \olim_m^{T(\A\otimes\B)} \bigl(\eta_{\A\otimes\B}(a_n\otimes b_m)\bigr)
% \end{equation}
% (cf. \cref{def:bisigmanormal}).

\begin{remark}\label{rem:monoidalmonad_cpu}
	The extension $\bar T = \bar U \bar F\colon \Calgcpu \to \Calgcpu$ is also a symmetric monoidal monad by the same reasoning of \cref{prop:monoidalmonad}, since the Kleisli category of $\bar T$ is given by Pedersen--Baire envelopes with $\sigma$-cpu maps, and this was shown to be a symmetric monoidal category in \cite[Remark 5.3.3]{fritz2025categories}.
\end{remark}

To be sufficiently self-contained, we now provide some details on the categorical results later applied to ensure that $\SigmaC$ and $\SigmaCcpu$ are symmetric monoidal categories.
Although the main reference~\cite{seal2013tensors} deals with monoidal monads and monoidal categories, here we restrict to \emph{symmetric} monoidal structures, as the additional generality never surfaces in our investigation.
\begin{notation} We adopt the following notation.
	\begin{itemize}
		\item The category of Eilenberg--Moore algebras of a monad $T\colon \clK\to \clK$ is denoted by $\EM(T)$ and its objects are denoted by pairs $(A,a)$ where $\A\in \clK$ and $a\colon T\A \to \A$ denotes the algebra structure.
		\item A symmetric monoidal category is a triple $(\clK, \otimes , I)$ where $\clK$ is a category, $\otimes \colon \clK \times \clK \to \clK$ denotes the tensor (monoidal) product, and $I$ denotes the unit: $\A \otimes I \cong \A$ for all $\A \in \clK$.
		\item A symmetric monoidal monad is a pair $(T,\theta)$ where $T$ is a monad on a symmetric monoidal category $(\clK,\otimes ,I)$, and $\theta$ is a natural transformation given by $\theta_{\A,\B} \colon T\A \otimes T\B \to T(\A\otimes \B)$.
	\end{itemize}
\end{notation}

\begin{definition}\label{def:kocktensor}
	Let $(T,\theta)$ be a symmetric monoidal monad on a symmetric monoidal category $(\clK,\otimes,I)$.
	The \newterm{tensor product} of Eilenberg--Moore algebras $(\A,a)$ and $(\B,b)$ is given by the coequalizer diagram
	\begin{equation*}%\label{eq:kocktensor}
		\xymatrix@C=5em{
		T(T\A\otimes T\B)
		\ar@<.7ex>[r]^{\mu T\theta }
		\ar@<-.7ex>[r]_{T(a\otimes b)}
		& T(\A\otimes \B)\ar[r]^{q_{(\A,a),(\B,b)}} & (\A,a)\kocktensor (\B,b)
		}
	\end{equation*}
	in $\EM(T)$.
	When the algebra morphisms are clear from context (or irrelevant), we simply write $\A \kocktensor \B$ and $q_{\A,\B}$.

	We say that $\EM(T)$ \newterm{has all tensors} if $q_{\A,\B}$ and $\A \kocktensor \B$ exist for all $(\A,a)$ and $(\B,b)$ in $\EM(T)$.
\end{definition}
\begin{remark}\label{rem:coeq_homtensor}
	\begin{itemize}
		\item  The coequalizer in the definition is reflexive with common section $T(\eta_A\otimes \eta_B)$.
		\item Although additional requests are needed to ensure one has a symmetric monoidal category (e.g., $\otimes$ should be associative), the above definition is sufficient to ensure the existence of $\phi\kocktensor \psi$ for morphisms $\phi\colon (\A,a) \to (\C,c)$ and $\psi\colon (\B,b)\to (\D,d)$, because of coequalizer properties:
		      \[
			      \xymatrix@C=5em{
			      T(\A\otimes \B)\ar[r]^{q_{\A,\B}}\ar[d]_{T(\phi\otimes \psi)} & \A\kocktensor \B\ar@{-->}[d]^{\phi\kocktensor \psi} \\
			      T(\C\otimes \D)\ar[r]^{q_{\C,\D}} & \C\kocktensor \D
			      }
		      \]
	\end{itemize}

\end{remark}

\begin{proposition}[{\cite[Lemma 2.3.3]{seal2013tensors}}]\label{prop:bijection_bimorphisms} Let $(T,\theta)$ be a symmetric monoidal monad on a symmetric monoidal category $(\clK,\otimes,I)$ such that $\EM(T)$ has all tensors.
	Consider $(\A,a)$, $(\B,b)$ and $(\C,c)$ Eilenberg--Moore algebras.
	The morphism
	\[
		\xymatrixcolsep{3pc}\xymatrix{
		\A \otimes \B \ar[r]^(0.45){\eta_{\A\otimes \B}}
		& T(\A\otimes \B) \ar[r]^-{q_{\A,\B}}
		& \A\kocktensor\B
		}
	\]
	induces, by precomposition, a bijection between the following:
	\begin{itemize}
		\item Morphisms of type $\A \kocktensor \B \to \C$.
		\item Bimorphisms of $\EM(T)$, i.e.\ morphisms $\omega\colon \A\otimes \B\to \C$ in $\clK$ satisfying
		      \[
			      \xymatrixcolsep{3pc}\xymatrix{
			      T\A \otimes T\B \ar[r]^{\theta_{\A,\B}} \ar[d]_{a\otimes b}
			      & T(\A\otimes \B) \ar[r]^-{T\omega}
			      & T\C \ar[d]^{c} \\
			      \A\otimes \B \ar[rr]_{\omega}
			      & & \C
			      }
		      \]
	\end{itemize}
\end{proposition}
\begin{theorem}[{\cite[Theorem 2.5.5]{seal2013tensors}}]\label{thm:SMC1}
	Let $(T,\theta)$ be a symmetric monoidal monad on a symmetric monoidal category $(\clK,\otimes,I)$ such that $\EM(T)$ has all tensors.
	If for all Eilenberg--Moore algebras $(\A,a)$, $(\B,b)$, $(\C,c)$, and $(\D,d)$,
	\begin{enumerate}
		\item \label{it:coeq_associators}
		      The diagrams of the form
		      \begin{equation*}
			      \begin{array}{c}
				      \xymatrix@C=6em{
				      T(T\A\otimes T\B)\kocktensor TT\C
				      \ar@<.7ex>[r]^{\mu T\theta \kocktensor \mu}
				      \ar@<-.7ex>[r]_{T(a\otimes b)\kocktensor Tc}
				       & T(\A\otimes \B)\kocktensor T\C \ar[r]^{q\kocktensor c}  & (\A\kocktensor \B)\kocktensor \C
				      }                                                                                             \\
				      \xymatrix@C=6em{
				      TT\A \kocktensor T(T\B\otimes T\C)
				      \ar@<.7ex>[r]^{\mu  \kocktensor \mu T\theta }
				      \ar@<-.7ex>[r]_{Ta \kocktensor T(b\otimes c)}
				       & T\A \kocktensor T(\B\otimes \C) \ar[r]^{a\kocktensor q} & \A\kocktensor (\B\kocktensor \C)
				      }
			      \end{array}
		      \end{equation*}
		      are coequalizers in $\EM(T)$;
		\item \label{it:coeq_tensor_epic}
		      $q_{\A,\B}\kocktensor q_{\C,\D}$ are epimorphisms in $\EM(T)$;
	\end{enumerate}
	Then $(\EM(T),\kocktensor,TI)$ is a symmetric monoidal category.
\end{theorem}

\begin{theorem}[{\cite[Theorem 2.6.4]{seal2013tensors}}]\label{thm:SMC2}
	Let $(T,\theta)$ be a symmetric monoidal monad on a symmetric monoidal category $(\clK,\otimes ,I)$ such that $\EM(T)$ has all tensors.
	If $T(\A\otimes -)$ preserves reflexive coequalizers in $\clK$ for all $\A \in \clK$, then $(\EM(T),\kocktensor,TI)$ is a symmetric monoidal category.
\end{theorem}
As it will be important in the proof of \cref{thm:tensor_scpu}, we mention that the proof of \cref{thm:SMC2} consists in showing that the assumptions of \cref{thm:SMC1} are satisfied.

\begin{proposition}[{\cite[Theorem 3.3]{khoshkam2002categorical}}]\label{prop:tensor_colims}
	Let $\A\in \Calg$.
	Maximal tensoring $\A\otimes - \colon \Calg\to \Calg$ preserves connected colimits.
\end{proposition}
%This fact seems somehow known, although it is hard to find an explicit proof.\footnote{One considering inductive limits can be found at \url{https://mathoverflow.net/questions/381936}.} 
For completeness, we provide a direct reasoning using the universal property of the maximal tensor product.
\begin{proof}
	Let us consider a connected category $\clC$ and a functor $L:\clC\to \Calg$ admitting a colimit $\colim L$.
	We now consider a cocone $\phi_c \colon Lc\otimes \A \to \B$ of $L(-)\otimes \A$ for some $\B\in \Calg$.
	We seek to show that there exists a compatible map $\colim L \otimes \A\to \B$, thereby ensuring that $\colim(L(-)\otimes A)\cong (\colim L)\otimes A$.

	Since $\clC$ is connected, one has a zig-zag of morphisms from any $c$ to $c'$.
	For simplicity, let us consider the case $f\colon c\to c'$.
	The commutative diagram
	\[
		\xymatrix{
		Lc\otimes \A \ar[rr]^{Lf \otimes \id_{\A}} \ar[dr]_{\phi_c}
		&&
		Lc'\otimes \A \ar[dl]^{\phi_{c'}}
		\\
		&
		\B
		&
		}
	\]
	ensures that $\phi_c(1\otimes a) = \phi_{c'}(1\otimes a)$.
	Since $\clC$ is connected, we obtain a unique map $\psi\colon \A \to \B$, which sends $a$ to $\phi_{c}(1\otimes a)$ for any choice of $c$.
	By looking at the marginalization on the other side, $\phi_c(-\otimes 1)\colon Lc\to \B$ induce a cocone for $L$. Therefore, one obtains $\phi\colon \colim L \to \B$.

	We now aim to show that the image of $\phi$ and $\psi$ commute. Once this is established, the universal property of the maximal tensor product would immediately yield $\colim L \otimes \A \to \B$, thus concluding the proof.

	Let us consider $x\in \colim L$. If $x$ is in the image of some $Lc$, then $\phi(x)$ commutes with the image of $\psi$ by construction.
	If not, $x$ is a uniform limit of a sequence of elements commuting with the image of $\psi$, because colimits in $\Calg$ are given by the completion of the colimits of $*$-algebras over $\bbC$ (see \cite[Corollaries 7.22 and 7.23]{bunke2021lectures}).\footnote{This fact was already used in \cref{lem:technical_access}.}
	In this case, one notes that $\phi(x)\psi(a)$ and $\psi(a)\phi(x)$ are both limits of $\phi(x_n)\psi(a)=\psi(a)\phi(x_n)$, and therefore they coincide.
	Finally, $\phi$ and $\psi$ commute, and the statement follows.
\end{proof}

\begin{theorem}\label{thm:tensor_shom}
	$\SigmaC$ is a symmetric monoidal category whose monoidal unit is $\bbC$ and whose tensor product $\kocktensor$ satisfies:
	\begin{enumerate}
		\item\label{it:free_algebras} $F(\A)\kocktensor F(\B)\cong F(\A\otimes \B)$;
		\item\label{it:universal_prop} The following universal property:
		      for every $\phi \colon \A \to \C$ and $\psi\colon \B \to \C$ in $\SigmaC$ with commuting images, there exists a unique $(\phi,\psi)\colon \A \kocktensor \B \to \C$ satisfying $(\phi,\psi)(a\kocktensor b) = \phi(a)\psi(b)$.
		      Moreover, all \starshom{}s of type $\A \kocktensor \B \to \C$ are obtained in this way.
	\end{enumerate}
\end{theorem}
\begin{proof}
	From \cref{prop:cocomplete,prop:tensor_colims}, $T(\A \otimes -)$ preserve reflexive coequalizers, and moreover $\SigmaC$ has tensors because it is cocomplete.
	By \cref{thm:SMC2}, $\SigmaC$ is a symmetric monoidal category (recall that $T$ is a symmetric monoidal monad by \cref{prop:monoidalmonad}).
	The monoidal unit is $\bbC$ because it is isomorphic to $T(\bbC)$ (for every finite-dimensional \cstar{} $\A$, $T(\A)\cong \A$, see \cite[Example 3.3.16]{fritz2025categories}).
	Item~\ref{it:free_algebras} follows from \cite[2.5.3]{seal2013tensors}.

	We are therefore reduced to show Item~\ref{it:universal_prop}.
	By \cref{prop:bijection_bimorphisms}, the morphisms $\A \kocktensor \B \to \C$ are in bijective correspondence with the bimorphisms $\A\otimes \B\to\C$ of $\EM(T)$ via precomposing with $q_{\A,\B}\eta_{\A\otimes \B}$.
	We claim that such morphisms exactly correspond to a choice of $\phi$ and $\psi$ as stated in Item~\ref{it:universal_prop}, which will then conclude the proof.

	Starting with a bimorphism $\omega\colon \A\otimes \B\to \C$, $\omega$ is uniquely determined by its marginals $\phi=\omega(-\otimes 1)$ and $\psi=\omega(1\otimes -)$ by the universal property of the maximal tensor product.
	Moreover, $\phi$ and $\psi$ have commuting images by direct check:
	\begin{equation*}
		\begin{split}
			\phi(a)\psi(b)&= \omega(a\otimes 1)\omega(1\otimes b) = \omega(a\otimes b) \\
			&= \omega((1\otimes b)(a\otimes 1))= \omega(1\otimes b)\omega(a\otimes 1) \\
			&= \psi(b)\phi(a).
		\end{split}
	\end{equation*}
	It then suffices to ensure that $\phi$ and $\psi$ so defined are necessarily $\sigma$-normal.
	For $\phi$, this follows by the commutative diagram
	\[
		\xymatrix@C=4em{
			\A\ar[r]^(0.45){\cong}\ar[rd]_{\cong}\ar@/^2em/[rrr]^{\phi} & \A \otimes \bbC \ar[r]^{\id_{\A}\otimes \, !}\ar[d]^{q\eta} & \A\otimes \B \ar[r]^-{\omega}\ar[d]^{q\eta}& \C\\
		&\A \kocktensor \bbC \ar[r]^{\id_{\A}\kocktensor \, !} & \A\kocktensor \B \ar[ru]&
		}
	\]
	since the maps in the lower path are \starshom{}s (where $!$ denotes the unique unital linear map from $\bbC$). Similarly, also $\psi$ is a \starshom.

	Conversely, let us consider two $\sigma$-homomorphisms $\phi$ and $\psi$ as in the statement.
	These yield a map $\omega \colon \A \otimes \B \to \C$ given by $\omega (a\otimes b)=\phi(a)\psi(b)$ via the universal property of the maximal tensor product. We claim that such an $\omega$ is a bimorphism in the sense of \cref{prop:bijection_bimorphisms}.
	By again employing the universal property of the maximal tensor product, the two compositions $T\A \otimes T\B\to \C$ can be tested on their decompositions $T\A \to \C$ and $T\B \to \C$.
	Since these are both $\sigma$-homomorphisms and $\C$ is a \cstars{}, both are uniquely determined by the $*$-homomorphisms $\A \to \C$ and $\B \to \C$ given by composing with the unit $\eta$.
	It is now a simple check that $\phi$ corresponds to both
	\[
		\A \to T(\A\otimes \B) \to T\C \to \C \qquad \text{and}\qquad \A \to \A\otimes \B\to \C
	\]
	and analogously for $\psi$, so the wanted commutation holds.
\end{proof}

\begin{proposition}\label{prop:tensor_commutative}
	Given $\A$, $\B$ commutative \cstars{}s, their tensor product $\A \kocktensor \B$ is again commutative.
	Moreover, $\A\kocktensor \B$ is the coproduct in $\cSigmaC$.
\end{proposition}

\begin{proof}
	Whenever $\A$ and $\B$ are commutative, $\A \otimes\B$ is commutative  (see, e.g., \cite[Lemma T.6.16]{weggeolsen1993Ktheory}).
	Moreover, $T$ preserves commutativity~\cite[Remark 5.3.4(iii)]{saito15monotone}, hence $T(\A\otimes \B)$ is also commutative.
	Since commutativity is preserved by quotients and $\A \kocktensor \B$ is a quotient of $T(\A\otimes \B)$ (by definition of the tensor and \cref{lem:reflexive_coeq}), we obtain that $\A\kocktensor \B$ is commutative.
	The second part of the statement is a direct application of the universal property of the tensor product (\cref{thm:tensor_shom}\ref{it:universal_prop}).
\end{proof}

\begin{remark}\label{rem:tensor_bool}
	Recall that $\cSigmaC$ is equivalent to the category of \bsigma{}s and $\sigma$-homomorphisms, as mentioned in \cref{sec:preliminaries}.
	By combining this equivalence with \cref{prop:tensor_commutative}, $\kocktensor$ retrieves the universal tensor product of \bsigma{}s \cite[Proposition 2.5.1]{fritz2025categories}.

	In particular, given two measurable spaces $X$ and $Y$, it is not always true that $\Linf(X)\kocktensor \Linf(Y)\cong \Linf(X\times Y)$ (cf.~\cite[Example 2.5.6]{fritz2025categories}, or the historical reference \cite[Example 37.A]{sikorski69boolean}).
	However, if $X$ and $Y$ are Baire measurable spaces (this notion was recalled before \cref{thm:giry_restriction}), we have that $\Linf(X)\cong F(C(X))$, and therefore
	\begin{align*}
		\Linf(X)\kocktensor \Linf(Y) &\cong F(C(X)) \kocktensor F(C(Y)) \\
		& \cong F(C(X)\otimes C(Y)) \\
		& \cong F(C(X\times Y)) \\
		& \cong \Linf(X\times Y)
	\end{align*}
	where $C(X)\otimes C(Y)\cong C(X\times Y)$ is well-known (see, e.g., \cite[Lemma T.6.16]{weggeolsen1993Ktheory}).
\end{remark}
% \begin{remark}
%   An alternative proof of \cref{prop:tensor_commutative}, suggested to us by Tobias Fritz, consists in noting that $\A^{\op}\kocktensor \B^{\op}\cong (\A\kocktensor \B)^{\op}$, where $-^{\op}$ denotes the opposite \cstar{} (i.e.\ where multiplication is reversed, $a\cdot^{\op} b\coloneqq a\cdot b$).
%   Indeed, the isomorphisms $\A^{\op} \otimes \B^{\op} \cong( \A\otimes \B)^{\op}$ and $T(\A^{\op})\cong T(\A)^{\op}$ are easily checked with sufficient grasp of the concepts. 
% \end{remark}

We now extend the tensor product to $\SigmaCcpu$.
As we will see, the proof boils down to noting that the coequalizers needed are reflexive and live in $\SigmaC$, and so \cref{lem:reflexive_coeq} applies.
\begin{theorem}\label{thm:tensor_scpu}
	$\SigmaCcpu$ is a symmetric monoidal category by extending the tensor product $\kocktensor$ of $\SigmaC$.
\end{theorem}
\begin{proof}
	We check that the hypotheses of \cref{thm:SMC1} are satisfied. Recall that $\bar{T}$ is a symmetric monoidal monad (\cref{rem:monoidalmonad_cpu}).

	We note that the algebra structures are uniquely determined, namely $a=U\epsilon_{\A}$, $b=U\epsilon_{\B}$ and $c=U\epsilon_{\C}$, as discussed at the end of the proof of \cref{thm:U_monadic}.
	In particular, the algebra structures are \starshom{}s, i.e.~they are in the image of the forgetful functor $\SigmaC\to \SigmaCcpu$.
	Since this functor preserves reflexive coequalizers by \cref{lem:reflexive_coeq}, we obtain that $\SigmaCcpu$ has tensors (the fact that tensors are reflexive coequalizer was noted in \cref{rem:coeq_homtensor}).

	To ensure that the diagrams of Item~\ref{it:coeq_associators} in \cref{thm:SMC1} are coequalizers, we first note that they already are so in $\SigmaC$ because \cref{thm:tensor_shom} applies \cref{thm:SMC2}, and the proof of this result consists in checking that the hypotheses of \cref{thm:SMC1} hold.
	Moreover, these coequalizer diagrams are reflexive: the common sections are $T(\eta_{\A}\otimes \eta_{\B})\kocktensor T(\eta_{\C})$ and $T(\eta_{\A})\kocktensor T(\eta_{\B}\otimes \eta_{\C})$, respectively.
	For the first coequalizer in Item~\ref{it:coeq_associators}, using functoriality of $T$ and $\kocktensor$, we have
		\begin{align*}
			(\mu T\theta \kocktensor \mu )(T(\eta_{\A}\otimes \eta_{\B})\kocktensor T(\eta_{\C}) ) & = \mu T(\theta (\eta_{\A}\otimes \eta_{\B}))\kocktensor \mu T(\eta_{\C})\\
			                                                                                       & = \mu T(\eta_{\A\otimes \B}) \kocktensor \id_{T\C}                      \\
			                                                                                       & = \id_{T(A\otimes B)}\kocktensor \id_{T\C}                              \\
			                                                                                       & = \id_{T(\A\otimes B)\kocktensor T\C}                                   \\
			(T(a\otimes b)\kocktensor Tc)(T(\eta_{\A}\otimes \eta_{\B})\kocktensor T(\eta_{\C}) )  & = T(a\eta_{\A}\otimes b\eta_{\B}) \kocktensor T(c\eta_{\C})             \\
			                                                                                       & = T(\id_{\A\otimes \B})\kocktensor T(\id_{\C})                          \\
			                                                                                       & = \id_{T(\A\otimes B)\kocktensor T\C}.
		\end{align*}
	An analogous reasoning applies to the second coequalizer in Item~\ref{it:coeq_associators}.
	As $\SigmaC \to \SigmaCcpu$ preserves reflexive coequalizers (\cref{lem:reflexive_coeq}), Item~\ref{it:coeq_associators} is verified.

	We now discuss Item~\ref{it:coeq_tensor_epic} of \cref{thm:SMC1}.
	We claim that $q_{\A,\B}\kocktensor q_{\C,\D}$ are in fact surjective maps.
	We first note that $q_{\A,\B}$ are all surjective by \cref{lem:reflexive_coeq}, and moreover there is a commutative diagram
	\[
		\xymatrix@C=5em{
		T(T(\A\otimes \B) \otimes T(\C\otimes \D))
		\ar[r]^{T(q\otimes q)}
		\ar@{->>}[d]_{q}
		&
		T(\A\kocktensor \B \otimes \C\kocktensor \D)
		\ar@{->>}[d]^{q}
		\\
		T(\A\otimes \B) \kocktensor T(\C\otimes \D)
		\ar[r]_{q\kocktensor q}
		&
		(\A\kocktensor \B) \kocktensor (\C\kocktensor \D)
		}
	\]
	In particular, the vertical arrows are surjective as already noted.
	To conclude, we only need to show that $T(q\otimes q)$ is also surjective. This follows from the fact that the maximal tensor product is exact~\cite[Proposition 3.7.1]{brownozawa} and $T$ preserves surjectivity (\cref{rem:T_surjective}).
\end{proof}

We now aim to describe all the \scpu{} maps of type $\A\kocktensor \B\to \C$.
To this end, we give the following.
\begin{definition}\label{def:bisigmanormal}
	Let us consider three \cstars{} $\A$, $\B$ and $\C$. A \cpu{} map $\omega\colon \A \otimes \B \to \C$ is \newterm{bi-$\sigma$-normal} if
	\[
		\omega\left( \olim_n^{\A} a_n \otimes \olim_m^{\B} b_m \right) = \olim_n^{\C} \olim_m^{\C} \omega(a_n \otimes b_m)
	\]
	for all order convergent sequences $(a_n)\subseteq \A$ and $(b_m)\subseteq \B$.
\end{definition}

\begin{lemma}\label{lem:theta_bisigma}
	The map $\theta\colon T\A \otimes T\B \to T(\A\otimes \B)$ is bi-$\sigma$-normal.
\end{lemma}
\begin{proof}
	By construction, $\theta$ is a $*$-homomorphism, so $\theta(a\otimes b)=\theta(a\otimes 1)\theta(1\otimes b)$, and $\theta(-\otimes 1)\colon T\A \to T(\A\otimes \B)$ is a $\sigma$-homomorphism, so $\theta(\olim_n a_n \otimes 1) = \olim_n \theta(a_n \otimes 1)$, and similarly for the second argument.

	Consequently, the statement holds by compatibility of order limits with multiplication:
	\begin{equation*}
		\begin{split}
			\theta\left(\olim_n^{T\A} a_n \otimes \olim^{T\B}_m b_m\right) & =\theta\left(\olim^{T\A}_n a_n \otimes 1\right)\theta\left(1\otimes \olim^{T\B}_m b_m\right)    \\
			                                                               & = \olim^{T(\A\otimes \B)}_n \theta(a_n \otimes 1)\olim^{T(\A\otimes \B)}_m \theta(1\otimes b_m) \\
			                                                               & = \olim^{T(\A\otimes \B)}_n \olim^{T(\A\otimes \B)}_m \theta(a_n \otimes 1)\theta(1\otimes b_m) \\
			                                                               & =\olim^{T(\A\otimes \B)}_n \olim^{T(\A\otimes \B)}_m \theta(a_n \otimes b_m).
		\end{split}
	\end{equation*}
\end{proof}
\begin{proposition}[Universal property of $\kocktensor$]\label{prop:universal_bisigma}
	The morphisms $\A \kocktensor \B \to \C$ in $\SigmaCcpu$ are in bijective correspondence to those $\omega\colon \A \otimes \B \to \C$ in $\Calgcpu$ that are bi-$\sigma$-normal.
\end{proposition}
\begin{proof}
	It suffices to show that $\omega\colon \A \otimes \B\to \C$ is bi-$\sigma$-normal if and only if is a bimorphism in the sense of \cref{prop:bijection_bimorphisms}.
	As already mentioned in the proof of \cref{thm:tensor_scpu}, the algebra structures are uniquely determined as shown at the end of the proof of \cref{thm:U_monadic}.
	We can therefore assume $a=U\epsilon_{\A}$, $b=U\epsilon_{\B}$ and $c=U\epsilon_{\C}$ in the bimorphism diagram of \cref{prop:bijection_bimorphisms}.

	If $\omega$ is a bimorphism, we take order convergent sequences $(a_n)\subseteq \A $ and $(b_m)\subseteq \B$.
	We wish to unravel the commuting diagram describing bimorphisms on the element $\olim_n \eta_{\A}(a_n)\otimes \olim_m \eta_{\B}(b_m)$.

	Concerning the upper path, we have
	\begin{align*}
		U\epsilon_{\C} T\omega \theta_{\A,\B} (\olim_n \eta_{\A}(a_n)\otimes \olim_m \eta_{\B}(b_m)) & = U\epsilon_{\C} T\omega (\olim_n\olim_m \theta_{\A,\B} (\eta_{\A}(a_n)\otimes \eta_{\B}(b_m)))  \\
		&= \olim_{n}\olim_m  U\epsilon_{\C} T\omega \theta_{\A,\B} (\eta_{\A}(a_n)\otimes \eta_{\B}(b_m))
	\end{align*}
	where we used \cref{lem:theta_bisigma} and $\sigma$-normality of $T\omega$ and $U\epsilon_{\C}$.
	Moreover, the following diagram commutes by unit compatibility of $\theta$ and naturality of $\eta$:
	\begin{equation}\label{eq:omega_comm}
		\xymatrixcolsep{3pc}\xymatrix{
		&\A\otimes \B\ar[ld]_{\eta_{\A}\otimes \eta_{\B}} \ar[d]^{\eta_{\A\otimes \B}} \ar[r]^{\omega}& \C\ar[d]^{\eta_{\C}}\ar[rd]^{\id_{\C}}\\
		T\A \otimes T\B \ar[r]^{\theta_{\A,\B}}
		& T(\A\otimes \B) \ar[r]^-{T\omega}
		& T\C \ar[r]^{U\epsilon_{\C}} & \C
		}
	\end{equation}
	Putting the two together, the upper path yields $\olim_n \olim_m \omega (a_n \otimes b_m)$.

	For the lower path,
	\[
		\omega (U(\epsilon_{\A})\otimes U(\epsilon_{\B} )(\olim_n \eta_{\A}(a_n)\otimes \olim_m \eta_{\B}(b_m)) = \omega (\olim_n a_n \otimes \olim_m b_m)
	\]
	by $\sigma$-normality of $U(\epsilon_{\A})$ and $U(\epsilon_{\B})$.
	Finally, by commutativity of the bimorphism diagram in \cref{prop:bijection_bimorphisms},
	\[
		\omega (\olim_n a_n \otimes \olim_m b_m) = \olim_n \olim_m \omega (a_n \otimes b_m),
	\]
	i.e.\ $\omega$ is bi-$\sigma$-normal.

	We now seek to prove the converse, i.e.\ that any bi-$\sigma$-normal ${\omega}$ is a bimorphism.
	This is done by two $\pi$-$\lambda$ arguments.
	For the first one, let
	\[
		M_{\B} \coloneqq \lbrace a \in T\A \mid \omega(U\epsilon_{\A}\otimes U\epsilon_{\B})(a \otimes \eta_{\B}(b)) = U\epsilon_{\C} T\omega \theta_{\A,\B}(a \otimes \eta_{\B}(b)) \text{ for all }b\in\B \rbrace\subseteq T\A
	\]
	This contains $\eta_{\A}(\A)$, because precomposing with $\eta_{\A}\otimes \id_{\B}$ yields $\omega$ on both sides of the defining equation by commutativity of \eqref{eq:omega_comm}.
	Moreover, $M_{\B}$ is a $*$-subspace because all maps involved are $*$-preserving and linear.
	Next, if $(a_n)\subseteq T\A$ is a bounded monotone increasing sequence contained in $M_{\B}$, then $\sup_n^{T\A} (a_n)$ belongs to $M_{\B}$ by bi-$\sigma$-normality of $\omega$ and $\theta$ (established in \cref{lem:theta_bisigma}).
	Since $T\A$ is the $\sigma$-closure of $\A$, we conclude that $M_{\B}=T\A$.

	For the second $\pi$-$\lambda$ argument, consider the $*$-subspace
	\[
		M_{T\A} \coloneqq \lbrace b \in T\B \mid \omega(U\epsilon_{\A}\otimes U\epsilon_{\B})(a \otimes b) = U\epsilon_{\C} T\omega\theta_{\A,\B}(a \otimes b) \text{ for all }a\in T\A\rbrace\subseteq T\B
	\]
	By the previous discussion, $\eta_{\B}(\B)\subseteq M_{T\A}$.
	Moreover, the argument above can be adapted to show that $M_{T\A}$ also is a $\sigma$-subspace, so $M_{T\A}=T\B$.
	Finally, the set
	\[
		M\coloneqq \lbrace z \in T\A\otimes T\B \mid \omega(U\epsilon_{\A}\otimes U\epsilon_{\B})(z) = U\epsilon_{\C} T\omega\theta_{\A,\B}(z)\rbrace\subseteq T\A\otimes T\B
	\]
	contains all simple tensors $z=a\otimes b$.
	$M$ is a norm-closed $*$-subspace, where norm-closure follows from the fact that all the maps are contractive.
	Therefore, $M$ contains all the norm limits of linear combinations of simple tensors, hence $M=T\A\otimes T\B$.
\end{proof}

% \begin{remark}
%   %This remark requires a hierarchical reasoning... Maybe not that important here, also the machinery would have to be introduced exactly for this remark, doesn't seem worth it.
%   A more explicit description of $\A \kocktensor \B$ is given by taking $T(\A\otimes \B)$ and quotienting it via the $\sigma$-ideal generated by
%   \[
%   \Im (\delta - \eta_{\A \otimes \B}(U\epsilon_{\A}\otimes U\epsilon_{\B}))
%   \]
%   This is an immediate consequence of the fact that the $\sigma$-ideal considered to define the tensor product is $\Im(\mu T\delta - T(U\epsilon_{\A}\otimes U\epsilon_{\B}))$, see \eqref{eq:tensor_def} in the proof of \cref{thm:tensor_scpu}.
%   Since $\mu T\delta $ and $T(U\epsilon_{\A}\otimes U\epsilon_{\B})$ are $\sigma$-homomorphisms, they are uniquely determined by their composition with the unit $\eta_{T\A\otimes T\B}$, from which one obtains
%   \[
%   \mu T(\delta) \eta_{T\A\otimes T\B} = \mu \eta_{T(\A\otimes \B)} \delta = \delta \qquad\text{ and }\qquad T(U\epsilon_{\A}\otimes U\epsilon_{\B})\eta_{T\A\otimes T\B} = \eta_{\A\otimes \B} (U\epsilon_{\A}\otimes U\epsilon_{\B}),
%   \] 
%   where both equalities follow from naturality of $\eta$.
% \end{remark}

%% file: sections/07_sigma_cp_instruments.tex
\section{\texorpdfstring{$\sigma$CP}{σCP}-instruments and tensor products}\label{sec:sigmaCP}
Following the standard formalism of von Neumann algebras, \cite{okamura2015measurement} discusses the notion of CP-instruments (\cite{davies1970operational,ozawa1984quantummeasuring}) and how these relate to tensor products of von Neumann algebras.
Here, we provide an analogue of their findings for \cstars{}s.

To fix the setting, recall that a CP-instrument assigns to each measurable subset of a measurable space $(X,\Sigma_X)$ a completely positive endomap of a chosen space, subject to a normalization condition and an adapted form of $\sigma$-additivity. 
In \cite{okamura2015measurement}, the chosen space may be taken to be a
von Neumann algebra $N$, and the endomaps to additionally be \emph{normal}.\footnote{In the reference, this is formalized in a slightly different fashion: the chosen space is the predual, where complete positivity still makes sense. However, since these maps bijectively correspond to normal completely positive maps on the von Neumann algebra, the present formulation is perfectly analogous.}

Here, we adapt this idea to \cstars{}s and \scpu{} maps, and obtain \cref{thm:sigmaCPinstruments}. 
In particular, we prove that $\sigmaCP$-instruments $\Sigma_X \to \sigmaCP(\A,\A)$ are in bijective correspondence with \scpu{} maps $\Linf(X)\kocktensor \A\to \A$. 
We emphasize that the considered maps are not only completely positive, but also $\sigma$-normal, and moreover the domain $\Linf(X)\kocktensor \A$ is itself a \cstars{}.

At first sight, our theorem is a \cstars{} counterpart of \cite[Proposition 3.3]{okamura2015measurement}, which states that every CP-instrument uniquely corresponds to a completely positive map out of a tensor product. 
However, the two are not fully parallel.
As the involved algebras in the reference are von Neumann, one might expect the map to be normal and the tensor product to be a von Neumann algebra, in analogy with \cref{thm:sigmaCPinstruments}.
Instead, the result yields a non-normal map on a
$C^*$-tensor product, so our statement is stronger.

These additional normality properties are not only important from our current perspective. 
In general, the completely positive map of \cite[Proposition 3.3]{okamura2015measurement} need not admit a normal extension to the von Neumann tensor product (\cite[\S V]{okamura2015measurement}), and this obstructs the construction of an associated measuring process (that is, the concrete dynamical realization of an instrument).
To overcome this, the authors introduced the \emph{normal extension property} for instruments, which requires precisely that such an extension exist.
Under this property, \cite[Theorem 3.4]{okamura2015measurement}
establishes a bijective correspondence between CP-instruments and
\emph{normal} completely positive maps from the relevant (von Neumann)
tensor product.

For this reason, \cref{thm:sigmaCPinstruments} captures the spirit of
\cite[Theorem 3.4]{okamura2015measurement}.
Notably, no additional hypothesis in the style of the normal extension
property is needed in our setting: every $\sigmaCP$-instrument
corresponds to a \scpu{} map.
Nonetheless, we make no claim about measuring processes. To the best
of our knowledge, relating our result to them would require a \cstars{}
version of the Stinespring's dilation theorem~\cite[Theorem 4.1]{paulsen2002}, which is not yet
available in the literature.

We expect the findings of this section to be of interest to current
developments in quantum measurement theory and instrument monads, as
they are closely related to the recent preprints
\cite{booth2026quantuminstruments,fritz2026quantuminstrumentmonad}.
In \cite{booth2026quantuminstruments}, the composition of CP-instruments relies on the normal extension property; in contrast,
\cite{fritz2026quantuminstrumentmonad} avoids these extension issues by
restricting to type I von Neumann algebras with separable preduals.
Although a conclusive comparison is beyond the present scope,
\cref{thm:sigmaCPinstruments} (together with
\cref{prop:instrument_composition}) suggests an alternative route in
which these restrictions are not needed, at the price of working with
\cstars{}s and $\sigma$-normal maps.

It is worth noting that, while the aforementioned references treat only measurable spaces, the present discussion applies to arbitrary \bsigma{}s.

%We are now ready to introduce the notion of $\sigmaCP$-instruments.
\begin{definition}
	Let $\Sigma$ be a \bsigma{} and let $\A $ and $\B$ be \cstars{}s. We denote with $\sigmaCP(\A,\B)$ the $\sigma$-normal completely positive maps $\A \to \B$.

	A \newterm{$\sigma$CP-measure} is a map $p\colon \Sigma \to \sigmaCP(\A,\B)$ such that:
	\begin{itemize}
		\item For all sequences $(U_n)_n\subseteq \Sigma$ of pairwise disjoint elements ($U_{n}\wedge U_{m}=\bot$ for all $n\neq m$) and all $a\in \A$, the sequence $(\sum_{n=1}^m p(U_n)(a))_m$ is order convergent in $\B$ and
		      \[\textstyle
			      p\left(\bigvee_n U_n\right)(a) = \sum_n p(U_n)(a)\coloneqq \olim_m^{\B} \left(\sum_{n=1}^m p(U_n)(a)\right)
		      \]
	\end{itemize}
	A \newterm{$\sigmaCP$-instrument} is a {$\sigma$CP-measure} such that $p(\top)(1)=1$.
\end{definition}
Alternatively, the property above may be reformulated by assuming $a\ge 0$ and using $\sup_n$ instead of $\olim_n$.

We now focus on preparatory material for the proof of the main theorem of this section (\cref{thm:sigmaCPinstruments}).

\begin{definition}
	Let $\Sigma$ be a \bsigma{} and let $\B$ be a \cstars{}. A (bounded) \newterm{$\B$-valued measure} is a map $q\colon \Sigma \to \B$ such that
	\begin{itemize}
		\item $q(U)\ge 0$ for all $U\in\Sigma$ and $q(\bot)=0$;
		\item For all sequences $(U_n)_n\subseteq \Sigma$ of pairwise disjoint elements, the sequence $(\sum_{n=1}^m q(U_n))_m$ is order convergent\footnote{Since $q(U_n)$ are positive and $\B$ is $\sigma$-closed, we can alternatively require the existence of $\lambda>0$ such that $\norm{\sum_{n=1}^m q(U_n)}\le \lambda$ for all $m$.} in $\B$
		      \[\textstyle
			      q\left(\bigvee_n U_n\right) = \sum_n q(U_n) \coloneqq \sup_m^{\B} \left(\sum_{n=1}^m q(U_n)\right)
		      \]
	\end{itemize}
\end{definition}
This definition is a generalization of \cite[Definition 2.1]{wright72measures} to the context of \bsigma{}s.
By additionally requiring $q(\top)=1$, one obtains a broader notion of positive operator-valued measures (see \cite[Example 3.2.8 and Definition 4.1.27]{fritz2025categories}).

A crucial remark is that for every \bsigma{} $\Sigma$, there is a valued measure $\Sigma\mapsto \Linfbool(\Sigma)$ denoted by $\chi_{-}$, which sends $U\in \Sigma$ to its abstract ``characteristic function''.
Concretely, $\chi_{U}\colon B(\bbC)\to \Sigma$ is given by $\chi_{U}(\lbrace 1\rbrace)=U$ and $\chi_{U}(\lbrace 0\rbrace)=\neg U$ (recall the description given in \cref{sec:preliminaries}).

In particular, for any sequence $(U_n)$,
\[
	\chi_{\bigvee_n U_n}=\sup_m\chi_{\bigvee_{n=1}^m U_n}\qquad\text{and}\qquad\chi_{\bigwedge_n U_n}= \inf_m \chi_{\bigwedge_{n=1}^m U_n}
\]
while for the bottom and top elements we have $\chi_{\bot}=0$ and $\chi_{\top}=1$.

For the categorically inclined, $\chi_{-}$ is the unit of the adjoint equivalence between $\cSigmaC$ and the category of \bsigma{}s with \bsigmahom{}s.
\begin{theorem}\label{thm:VM_sigmacp}
	Let $\Sigma$ be a \bsigma{} and $\B$ a \cstars{}.
	Then every $\B$-valued measure $q$ uniquely extends to a $\sigma$-normal completely positive map $\tilde{q}\colon \Linfbool(\Sigma)\to \B$ in the sense of the following commutative diagram:
	\[
		\xymatrix{\Sigma\ar[rr]^{\chi_{-}}\ar[rd]_{q}&&\Linfbool(\Sigma)\ar[ld]^{\widetilde{q}}\\
		&\B&
		}
	\]
\end{theorem}
\begin{proof}[Sketch of the proof]
	The argument is a simple rewriting of the proof of \cite[Theorem 4.1.30]{fritz2025categories}.
	Let us offer some details for the sake of clarity.
	If $\Sigma$ is an honest $\sigma$-algebra, i.e.\ there exists a measurable space $(X,\Sigma_X)$ with $\Sigma_X\cong \Sigma$, then this result was shown in \cite[Section 2]{wright72measures} and basically follows by rewriting the construction of the Lebesgue integral.
	For a general \bsigma{} $\Sigma$, by Loomis--Sikorski representation theorem \cite[Theorem 29.1]{sikorski69boolean}, there exists a measurable space $(X,\Sigma_X)$ such that $\Sigma\cong \Sigma_X/I$, with $I$ a $\sigma$-ideal.
	In particular, $\Linfbool(\Sigma)\cong \Linf(X)/\mathcal{I}$, where $\mathcal{I}$ is the set of bounded measurable functions whose support lies in $I$ (\cite[Proposition 4.1.22]{fritz2025categories}).
	We now consider
	\[
		\xymatrix{\Sigma_{X}
		\ar[rrrr]^{\chi_-}
		\ar@{->>}[dr]_{r}
		&&&&
		\Linf(X)
		\ar@{->>}[dl]^{\widetilde{\chi_{-} r}}
		\ar@/^2.5pc/[ddll]^{\widetilde{qr}}
		\\
		&
		\Sigma
		\ar[rr]^{\chi_-}
		\ar[dr]_{q}
		&&
		\Linf(X)/\mathcal{I}
		\ar@{-->}[dl]^{\widetilde{q}}
		\\
		&&
		\A}
	\]
	where $\widetilde{qr}$ exists from the previous step.
	We claim that the dashed map exists by the quotient universal property.
	Indeed, for any self-adjoint $g\in \mathcal{I}$, we have
	\[
		-\norm{g}\chi_{\operatorname{supp}(g)} \le g\le \norm{g} \chi_{\operatorname{supp}(g)}
	\] with $\operatorname{supp}(g)\in I$ the support of $g$ (multiply the standard inequality $-\norm{g}\le g\le \norm{g}$ by $\chi_{\operatorname{supp}(g)}$).
	Together with $\widetilde{qr}(\chi_{\operatorname{supp}(g)})=qr({\operatorname{supp}(g)})= 0$, this implies
	\[ 0 =-\norm{g}\cdot 0 = \widetilde{qr}(-\norm{g}\chi_{\operatorname{supp}(g)}) \le \widetilde{qr}(g)\le \widetilde{qr}(\norm{g} \chi_{\operatorname{supp}(g)}) = \norm{g} \cdot 0 =0,
	\]
	and therefore $\widetilde{qr}(g)=0$.
\end{proof}
For the next preliminary result, we recall the following.
\begin{definition}
	In a \cstar{} $\A$, a \newterm{projection} $p$ is an element of $\A$ such that $p^*=p=p^2$.

	A \newterm{simple element} $s$ is a linear combination of projections: $s=\sum_{k=1}^m \lambda_k p_k$.
\end{definition}
In $\Linf(X)$, a simple element is given by a simple function because projections are characteristic functions of measurable subsets.
Similarly, in $\Linfbool(\Sigma)$, projections are in correspondence with the elements of $\Sigma$ via $\chi_{-}\colon \Sigma \to \Linfbool(\Sigma)$ (\cite[Lemma 4.1.26]{fritz2025categories}).

\begin{lemma}\label{lem:technical_simple_els}
	Let $\Sigma$ be a \bsigma{} and let $\A$ be a \cstars{}.
	\begin{enumerate}
		\item\label{it:simple_els_dense} The linear span of elements $s\otimes a$, where $a\in \A$ and $s$ is a simple element in $\Linfbool(\Sigma)$, is dense in $\Linfbool(\Sigma)\otimes \A$.
		      Moreover, the self-adjoint elements of this span is dense in the set of self-adjoint elements of $\Linfbool(\Sigma)\otimes \A$.
		\item\label{it:simple_els_positivity} $\sum_{k=1}^m s_k \otimes a_k$ is positive if and only there exist some $n$ and some $t_{\ell}\ge 0$ and $b_{\ell}\ge 0$ for every $\ell \in \lbrace 1, \dots , n\rbrace$, such that
		      \[
			      \sum_{k=1}^m s_k \otimes a_k = \sum_{\ell=1}^n t_{\ell} \otimes b_{\ell}.
		      \]
		      Moreover, $t_{\ell}$ can be taken to be disjoint projections: $t_{\ell}\cdot t_{\ell'}=0$ for all distinct $\ell,\ell'$.
	\end{enumerate}
\end{lemma}
\begin{proof}
	Item~\ref{it:simple_els_dense} follows from the fact that linear combinations of $f\otimes a$ are dense by definition of the tensor product (since it is the norm-closure of the algebraic tensor product), and the fact that simple elements are dense in $\Linfbool(\Sigma)$ by \cite[Corollary 3.4.4]{fritz2025categories}.
	The restriction to self-adjoint elements is analogous.

	Regarding Item~\ref{it:simple_els_positivity}, let us assume that $\sum_{k=1}^m s_k \otimes a_k\ge 0$.
	Let $s_k= \sum_{j} \lambda_{j,k} p_{j,k}$ with $\lambda_{j,k}\in\bbC$ and $p_{j,k}$ projections.
	By standard properties of the algebraic tensor product,
	\[\textstyle
		\sum_k s_k \otimes a_k =  \sum_k \Big(\sum_j  \lambda_{j,k} p_{j,k}\Big)\otimes a_k = \sum_{j,k} p_{j,k}\otimes (\lambda_{j,k}a_k)
	\]
	Hence, we can assume without loss of generality that our element is written as $\sum_{k=1}^m p_k \otimes a_k$ with $p_k$ projections.
	For the next step, we introduce some auxiliary notation. Let $p_k^{\top}\coloneqq p_k$ and $p_k^{\bot}\coloneqq 1-p_k$.
	Now, for any Boolean sequence $I=(i_1, \dots, i_m) \in \lbrace \bot,\top\rbrace^m$, we define
	\[
		p^I \coloneqq (p_1)^{i_1} \cdot (p_2)^{i_2} \cdots (p_m)^{i_m}, \qquad \text{and}\qquad a^I \coloneqq \sum_{j \,\mid\, i_j=\top} a_j
	\]
	By a simple induction reasoning,
	\[
		\sum_{k=1}^m p_k \otimes a_k = \sum_{I \in \lbrace \bot,\top\rbrace^m} p^I \otimes a^I,
	\]
	and moreover $p^I \cdot p^J=0$ for any distinct $I,J$ because $p_k^{\top}p_k^{\bot} =0$ for any $k$.
	By restricting to those $I\in \lbrace \bot,\top\rbrace^m$ such that $p^I\neq 0$, we can set $t_{\ell}=p^I$ and $a^I=b^{\ell}$ and obtain the formula described in the statement.
	However, we have yet to show that $a^I$ is positive (while positivity of $p^I$ is immediate since it is a projection).
	To this end, we recall the well-known fact that, for any compact Hausdorff space $Y$, $C(Y)\otimes \A \cong C(Y,\A)$, where $C(Y,\A)$ is the \cstar{} given by continuous maps $Y\to \A$ where operations are defined pointwise (see, e.g., \cite[Corollary T.6.17]{weggeolsen1993Ktheory}).
	From this description, one notes that $f\in C(Y,\A)$ is positive if and only if $f(y)\ge 0$ for any $y\in Y$.
	By taking $Y$ the Gelfand spectrum of $\Linfbool(\Sigma)$ and $U_I $ the subset of $Y$ satisfying $p^I = \chi_{U_I}$,\footnote{$p^I$ interpreted as a function $Y\to\bbC$ has image contained in $\lbrace 0,1\rbrace$ because it is a projection.} for any element $u\in U^I$,
	\[
		p^J(u)=p^J(u)p^I(u)= 0
	\] for any $J\neq I$, because $p^J p^I=0$.
	Moreover, $U^I$ is nonempty because we restricted to those $I$ such that $p^I\neq 0$.
	It follows that $a^I$ is exactly the image of $u$ via the continuous map $\sum_{I \in \lbrace \bot,\top\rbrace^m} a^I \cdot p^I(-) \colon Y \to \A$, which is positive by assumption.
	Finally, $a^I$ is positive, and the statement is shown.
\end{proof}

\begin{theorem}\label{thm:sigmaCPinstruments}
	Let $\Sigma$ be a \bsigma{} and let $\A$ and $\B$ be \cstars{}s.
	There is a bijective correspondence between
	\begin{enumerate}
		\item\label{it:sigmaCPinstrument} $\sigmaCP$-instruments $\Sigma \to \sigmaCP(\A,\B)$;
		\item\label{it:scpu_tensor} \scpu{} maps $\Linfbool(\Sigma)\kocktensor \A \to \B$.
	\end{enumerate}
\end{theorem}
\begin{proof}
	We define a map $p_{-}$ from Item~\ref{it:scpu_tensor} to Item~\ref{it:sigmaCPinstrument} by sending a \scpu{} map $\phi\colon \Linfbool(\Sigma)\kocktensor \A \to \B$ to the $\sigmaCP$-instrument $p_{\phi}\colon U \mapsto (a\mapsto \phi(\chi_U \otimes a))$.
	This is indeed a $\sigmaCP$-measure: for every sequence $(U_n)\subseteq \Sigma$ of pairwise disjoint elements,
	\begin{align*}
		\textstyle p_{\phi}\left(\bigvee_n U_n\right)(a) & \textstyle =\phi(\chi_{\bigvee_n U_n} \otimes a)                                                                                                  \\
		                                                 & \textstyle = \phi \left(\left(\sup_m \sum_{n=1}^m\chi_{U_n} \right)\otimes a\right)   \tag{$\chi_{\bigvee_n U_n}=\sup_m \sum_{n=1}^m \chi_{U_n}$} \\
		                                                 & \textstyle = \olim_m \phi\left(\left(\sum_{n=1}^m\chi_{U_n} \right) \otimes a \right) \tag{\cref{prop:universal_bisigma}}                         \\
		                                                 & \textstyle = \olim_m \sum_{n=1}^m \phi(\chi_{U_n}\otimes a)                           \tag{linearity}                                             \\
		                                                 & \textstyle = \sum_{n} p_{\phi} (U_n)(a)
	\end{align*}
	Normalization is obvious: $p_{\phi}(\top)(1)=\phi(\chi_{\top}\otimes 1)= \phi(1\otimes 1)=1$.

	Injectivity of $p_{-}$ is also straightforward.
	If $p_{\phi}=p_{\psi}$, then $\phi(\chi_U\otimes a)=\psi(\chi_U\otimes a)$ for all $a\in \A$ and all $U\in\Sigma$. Since for every $f\in \Linfbool(\Sigma)$ one has $f=\olim_n s_n$, where $s_n$ are simple elements (\cite[Corollary 3.4.4]{fritz2025categories}), this is sufficient to conclude that $\phi$ and $\psi$ are equal on all simple tensors by bi-$\sigma$-normality (\cref{prop:universal_bisigma}), and therefore they must coincide.

	We are left to show surjectivity of $p_{-}$.
	We first note that, for any positive $a\in \A$, a $\sigmaCP$-instrument $p$ yields a $\B$-valued measure $q_a(U)\coloneqq p(U)(a)$.
	The $\sigma$-additivity of $q_a$ is immediate from that of $p$, so we only need to check $q_a(U)\ge 0$ and $q_a(\bot)=0$.
	The first one follows from the fact that $a$ is positive and $p(U)$ is completely positive.
	The second one is a consequence of $q_a(U\vee \bot)=q_a(U)+q_a(\bot)$ for all $U\in \Sigma$ by $\sigma$-additivity of $q_a$: indeed, this entails $q_a(\bot)=2q_a(\bot)$, and so $q_a(\bot)=0$.
	By \cref{thm:VM_sigmacp}, we obtain a $\sigma$-normal cp map $\Linfbool(\Sigma)\to \B$ extending $q_a$.
	Moreover, whenever $f=\olim_n s_n$ with $s_n=\sum_\ell \lambda_{\ell,n} \chi_{U_{\ell,n}}$ simple elements, one has $\widetilde{q_a} (f)= \olim_n \sum_\ell \lambda_{\ell,n} p({U_{\ell,n}})(a)$ by $\sigma$-normality of $\widetilde{q_a}$.
	Overall, we have obtained a map
	\[
		\phi \colon   \xymatrix@R=.1pt{
		\Linfbool(\Sigma)\times \A_+ \ar[r]& \B \\
		(f,a)\ar@{|->}[r]
		& \widetilde{q_a}(f)
		}
	\]
	where $\A_+$ denotes the positive elements of $\A$.
	By linearity, we can extend the whole discussion to all elements of $\A$.
	Explicitly, given $a=\sum_i \mu_i a_i = \sum_j \nu_j a_j'$ finite sums with $a_i,a_j'\ge 0$ and $\mu_i,\nu_j \in \bbC$, we have
	\begin{align*}
		\textstyle\sum_i\mu_i \phi(f,a_i) & =\textstyle \sum_i \mu_i \left(\olim_n \sum_{\ell} \lambda_{\ell,n} p(U_{\ell,n})(a_i)\right)          \\
		                                 & =\textstyle \olim_n \sum_{\ell} \lambda_{\ell,n}p(U_{\ell,n}) \left(\sum_i \mu_i a_i\right) \tag{\dag} \\
		                                 & =\textstyle\olim_n \sum_{\ell} \lambda_{\ell,n}p(U_{\ell,n}) \left(\sum_j \nu_j a_j'\right)            \\
		                                 & =\textstyle \sum_j \nu_j \olim_n  \sum_{\ell} \lambda_{\ell,n} p(U_{\ell,n})(a_j') \tag{\ddag}          \\
		                                 & =\textstyle\sum_j \nu_j \phi(f,a_j')
	\end{align*}
	where the two equalities $(\dag),(\ddag)$ hold by linearity of order limits and $p(U)$, for all $U \in \Sigma$.
	% In particular,
	% \[
	% \phi(f,a) = \olim_n \sum_{\ell} \lambda_{\ell,n} p(U_{\ell,n})(a)
	% \]
	% holds for all elements $a\in \A$ and all choices of $f= \olim_n \sum_\ell \lambda_{\ell,n} \chi_{U_{\ell,n}}$.
	Finally, we obtain a linear map $\Linfbool(\Sigma)\odot \A \to \B $, which we again call $\phi$ for the sake of simplicity, where $\odot$ denotes the algebraic tensor product.

	By \cref{lem:technical_simple_els}\ref{it:simple_els_positivity}, we see that this map restricted to the elements $\sum_{k=1}^m s_k\otimes a_k$ is positive, and therefore, by \cite[Proposition 2.1]{paulsen2002} and its proof, it is bounded and satisfies $\norm{\phi(b)}\le \norm{b}$ for all self-adjoint $b$.
	Finally, $\phi$ can be uniquely extended to the norm-completion of its domain, which coincides with $\Linfbool(\Sigma)\otimes \A$ by \cref{lem:technical_simple_els}\ref{it:simple_els_dense}.
	Moreover, such an extension must be self-adjoint, and contractive on self-adjoint elements, by the second sentence of \cref{lem:technical_simple_els}\ref{it:simple_els_dense}.
	We can now apply \cite[Corollary 3.25]{fritz2023involutive} to conclude that this map is positive.

	By reapplying the same strategy to the $\sigmaCP$-instrument $\Sigma \to \sigmaCP(\A\otimes M_n , \B\otimes M_n)$ sending $U$ to $p(U)\otimes \id_{M_n}$, we obtain positivity of $\Linfbool(\Sigma)\otimes \A\otimes M_n \to \B \otimes M_n$.
	Arbitrariness of $n$ ensures complete positivity.

	We now check bi-$\sigma$-normality of $\phi\colon \Linfbool(\Sigma)\otimes \A\to \B$.
	On $\Linfbool(\Sigma)$, $\sigma$-normality follows from construction: given $a\in \A$ and $a=\sum_i \mu_i a_i$ with $a_i\ge 0$ and $\mu_i \in \bbC$,
	\begin{align*}
		\phi(\olim_k f_k \otimes a) & =\sum_i \mu_i \phi(\olim_k f_k \otimes a_i) = \sum_i \mu_i q_{a_i}(\olim_k f_k)  \\
		                            & = \sum_i \mu_i \olim_k q_{a_i}(f_k) = \olim_k \sum_i \mu_i \phi(f_k \otimes a_i) \\
		                            & = \olim_k \phi(f_k\otimes a)
	\end{align*}
	(see also \cite[Theorem 4.1.7]{fritz2025categories}).

	The other side is more delicate.
	Let us restrict to the case of $f\in \Linfbool(\Sigma)$ positive, and let $(a_n)\subseteq \A $ be a monotone decreasing sequence such that $\inf_n a_n =0$; in particular, $a_n\ge 0$ for all $n$.
	Then, $\phi(f\otimes a_n)\le \phi(f\otimes a_{n-1})$, and we wish to show that $\inf_n \phi(f\otimes a_n) =0$.
	To this end, note that $0 \le \phi(f\otimes a_n)\le\phi(\norm{f}\otimes a_n)$ by positivity of $\phi$ and the standard inequality $0\le f\le \norm{f}$.
	Consequently,
	\[\textstyle
		0 \le \inf_n \phi(f\otimes a_n) \le \inf_n \phi(\norm{f}\otimes a_n) \overset{(\dag)}{=} \inf_n \norm{f} p(\top)(a_n) = \norm{f}p(\top)(\inf_n a_n) = 0
	\]
	where the equality $(\dag)$ follows from definition of $\phi$, while the next equality holds by $\sigma$-normality of $p(\top)$. Therefore, $\inf_n \phi(f\otimes a_n)=0=\phi(f\otimes \inf_n a_n)$, as wanted. By linearity, $\phi(f\otimes -)$ also preserves any suprema of norm-bounded monotone increasing sequences (cf.\ \cite[Remark 3.2.5]{fritz2025categories}), and so it preserves also any order limit (as mentioned in \cref{sec:preliminaries}).
	Since every element of a \cstar{} is a linear combination of positive elements, $\phi(f \otimes \olim_m a_m) = \olim_m \phi(f\otimes a_m)$ holds for all $f\in \Linfbool(\Sigma)$.

	Finally, bi-$\sigma$-normality of $\phi$ is ensured, and surjectivity of $p_{-}$ now follows from \cref{prop:universal_bisigma}.
\end{proof}

We close with the following composition of $\sigmaCP$-instruments,
a first step toward a variant for \cstars{}s of the program of quantum
instrument monads~\cite{booth2026quantuminstruments,fritz2026quantuminstrumentmonad}.

\begin{proposition}\label{prop:instrument_composition}
		Let $\Sigma$ and $\Sigma'$ be \bsigma{}s and let $\A$, $\B$, $\C$ be \cstars{}s.
		$\sigmaCP$-instruments $p\colon \Sigma \to \sigmaCP(\A,\B)$ and $p'\colon \Sigma' \to \sigmaCP(\B,\C)$ compose to a $\sigmaCP$-instrument
		\[
			p' \circ p \colon \Sigma \otimes \Sigma' \longrightarrow \sigmaCP(\A,\C)
		\]
		whose outcomes are indexed by the universal tensor product of \bsigma{}s.
\end{proposition}
\begin{proof}
		By \cref{thm:sigmaCPinstruments}, $p$ and $p'$ correspond to \scpu{} maps $\phi\colon \Linfbool(\Sigma)\kocktensor \A \to \B$ and $\phi'\colon \Linfbool(\Sigma')\kocktensor \B \to \C$.
		Since $\SigmaCcpu$ is symmetric monoidal (\cref{thm:tensor_scpu}), the composite
		\[
			\Linfbool(\Sigma')\kocktensor \Linfbool(\Sigma)\kocktensor \A
			\xrightarrow{\ \id \kocktensor \phi\ }
			\Linfbool(\Sigma')\kocktensor \B
			\xrightarrow{\ \phi'\ }
			\C
		\]
		is again a \scpu{} map.
		By \cref{prop:tensor_commutative} and \cref{rem:tensor_bool}, $\Linfbool(\Sigma')\kocktensor \Linfbool(\Sigma)\cong \Linfbool(\Sigma\otimes \Sigma')$, where $\otimes$ denotes the universal tensor product of \bsigma{}s~\cite[Proposition 2.5.1]{fritz2025categories}.
		Applying \cref{thm:sigmaCPinstruments} in the other direction yields the desired instrument.
\end{proof}

%% file: references.bib
@article{wright76regular,
author = {Wright, J. D. Maitland},
title = {Regular {$\sigma$}-Completions of {C}*-Algebras},
journal = {Journal of the London Mathematical Society},
volume = {s2-12},
number = {3},
pages = {299-309},
doi = {https://doi.org/10.1112/jlms/s2-12.3.299},
year = {1976}
}

@book{saito15monotone,
  publisher = {Springer London},
  year = {2015},
  Author = {Sait{\^o}, Kazuyuki and Wright, J. D. Maitland},
 Title = {Monotone {C}omplete $C^*$-algebras and Generic Dynamics},
 Series = {Springer Monogr. Math.},
 DOI = {https://doi.org/10.1007/978-1-4471-6775-4}
}

@book {sikorski69boolean,
    AUTHOR = {Sikorski, Roman},
     TITLE = {Boolean algebras},
    SERIES = {Ergebnisse der Mathematik und ihrer Grenzgebiete [Results in
              Mathematics and Related Areas]},
    VOLUME = {25},
   EDITION = {Third},
 PUBLISHER = {Springer-Verlag New York, Inc., New York},
      YEAR = {1969},
     PAGES = {x+237}
}

@book {pedersen2018automgroups,
    AUTHOR = {Pedersen, Gert K.},
     TITLE = {{$C^*$}-algebras and their automorphism groups},
    SERIES = {Pure and Applied Mathematics (Amsterdam)},
   EDITION = {Second},
      NOTE = {Edited and with a preface by S\o ren Eilers and Dorte Olesen},
 PUBLISHER = {Academic Press, London},
      YEAR = {2018},
     PAGES = {xviii+520},
      ISBN = {978-0-12-814122-9}
}

@book {borceux94handbook2,
    AUTHOR = {Borceux, Francis},
     TITLE = {Handbook of categorical algebra. 2},
    SERIES = {Encyclopedia of Mathematics and its Applications},
    VOLUME = {51},
      NOTE = {Categories and structures},
 PUBLISHER = {Cambridge University Press, Cambridge},
      YEAR = {1994},
     PAGES = {xviii+443},
      ISBN = {0-521-44179-X},
   MRCLASS = {18-02 (18Exx)},
  MRNUMBER = {1313497},
MRREVIEWER = {Martin\ Hyland},
}

@book{johnstone82stone,
    AUTHOR = {Johnstone, Peter T.},
     TITLE = {Stone spaces},
    SERIES = {Cambridge Studies in Advanced Mathematics},
    VOLUME = {3},
 PUBLISHER = {Cambridge University Press, Cambridge},
      YEAR = {1982},
     PAGES = {xxi+370},
      ISBN = {0-521-23893-5},
   MRCLASS = {54-02 (06Dxx 18B30)},
  MRNUMBER = {698074},
MRREVIEWER = {Andreas\ Blass},
}

@InProceedings{giry82probability,
author="Giry, Mich{\`e}le",
editor="Banaschewski, B.",
title="A categorical approach to probability theory",
booktitle="Categorical Aspects of Topology and Analysis",
year="1982",
publisher="Springer Berlin Heidelberg",
address="Berlin, Heidelberg",
pages="68--85",
isbn="978-3-540-39041-1",
doi={10.1007/BFb0092872},
}

@Book{brownozawa,
 Author = {Brown, Nathanial P. and Ozawa, Narutaka},
 Title = {$C^*$-algebras and finite-dimensional approximations},
 FSeries = {Graduate Studies in Mathematics},
 Series = {Grad. Stud. Math.},
 ISSN = {1065-7338},
 Volume = {88},
 ISBN = {978-0-8218-4381-9},
 Year = {2008},
 Publisher = {Providence, RI: American Mathematical Society (AMS)},
 Language = {English}
}

@book{murphy90cstar,
	author = {Murphy, Gerard J.},
	title = {C*-algebras and operator theory},
	publisher = {Academic Press},
	year = {1990},
}

@article{wright72measures,
    author = {Wright, J. D. Maitland},
    title = "{Measures with Values in a Partially Ordered Vector Space}",
    journal = {Proceedings of the London Mathematical Society},
    volume = {s3-25},
    number = {4},
    pages = {675-688},
    year = {1972},
    month = {11},
    issn = {0024-6115},
    doi = {https://doi.org/10.1112/plms/s3-25.4.675},
    eprint = {https://academic.oup.com/plms/article-pdf/s3-25/4/675/4505769/s3-25-4-675.pdf},
}

@Book{paulsen2002,
 Author = {Paulsen, Vern},
 Title = {Completely bounded maps and operator algebras},
 FSeries = {Cambridge Studies in Advanced Mathematics},
 Series = {Camb. Stud. Adv. Math.},
 Volume = {78},
 ISBN = {0-521-81669-6},
 Year = {2002},
 Publisher = {Cambridge: Cambridge University Press},
 Language = {English},
}

@misc{fritz2023involutive,
  title={Involutive {M}arkov categories and the quantum de {F}inetti theorem},
  author={Fritz, Tobias and Lorenzin, Antonio},
	note = {\href{https://arxiv.org/abs/2312.09666}{arXiv:2312.09666}},
  year={2023}
}

@Article{jamneshan23foundational,
 Author = {Jamneshan, Asgar and Tao, Terence},
 Title = {Foundational aspects of uncountable measure theory: {Gelfand} duality, {Riesz} representation, canonical models, and canonical disintegration},
 FJournal = {Fundamenta Mathematicae},
 Journal = {Fundam. Math.},
 ISSN = {0016-2736},
 Volume = {261},
 Number = {1},
 Pages = {1--98},
 Year = {2023},
 Language = {English},
 DOI = {https://doi.org/10.4064/fm226-7-2022},
 note = {\href{https://arxiv.org/abs/2010.00681}{arXiv:2010.00681}},
}

@article{furber_jacobs_gelfand,
    author = {Furber, Robert W. J. and Jacobs, Bart P. F.},
    title  = {From {K}leisli Categories to Commutative {C}*-algebras: Probabilistic {G}elfand Duality},
    journal = {Logical Methods in Computer Science},
    volume = {11},
    issue = {2},
    year = {2015},
    doi = {10.2168/LMCS-11(2:5)2015},
    note = {\href{https://arxiv.org/abs/1303.1115}{arXiv:1303.1115}},
}

@incollection {milne_notes,
	AUTHOR = {Deligne, P.},
	TITLE = {Cat\'{e}gories tannakiennes},
	BOOKTITLE = {The {G}rothendieck {F}estschrift, {V}ol. {II}},
	SERIES = {Progr. Math.},
	VOLUME = {87},
	PAGES = {111--195},
	PUBLISHER = {Birkh\"{a}user Boston, Boston, MA},
	YEAR = {1990},
	MRCLASS = {14A99 (12H05 18A99)},
	MRNUMBER = {1106898},
	MRREVIEWER = {James Milne},
}

@article{Kock1970,
  title = {Monads on symmetric monoidal closed categories},
  volume = {21},
  ISSN = {1420-8938},
  url = {http://dx.doi.org/10.1007/BF01220868},
  DOI = {10.1007/bf01220868},
  number = {1},
  journal = {Archiv der Mathematik},
  publisher = {Springer Science and Business Media LLC},
  author = {Kock,  Anders},
  year = {1970},
  month = Dec,
  pages = {1–10}
}

@book {takesaki1979operatoralgebras1,
    AUTHOR = {Takesaki, Masamichi},
     TITLE = {Theory of operator algebras. {I}},
 PUBLISHER = {Springer-Verlag, New York-Heidelberg},
      YEAR = {1979},
     PAGES = {vii+415},
      ISBN = {0-387-90391-7},
   MRCLASS = {46Lxx},
  MRNUMBER = {548728},
MRREVIEWER = {Christopher\ Lance},
}

@misc {furber2024monad,
	author = {Furber, Robert},
	title = {A Probability Monad on Measure Spaces},
	note = {\href{https://www.robertfurber.com/preprints/cemonad-preprint.pdf}{robertfurber.com/preprints/cemonad-preprint.pdf}},
}

@article {pavlov2022duality,
	author = {Pavlov, Dmitri},
	title = {Gelfand-type duality for commutative von {N}eumann algebras},
	journal = {J. Pure Appl. Algebra},
	volume = {226},
	number = {4},
	pages = {106884},
	year = {2022},
	note = {\href{https://arxiv.org/abs/2005.05284}{arXiv:2005.05284}},
}

@book {EGNO,
	AUTHOR = {Etingof, P. and Gelaki, S. and Nikshych, D. and
	Ostrik, V.},
	TITLE = {Tensor categories},
	SERIES = {Mathematical Surveys and Monographs},
	VOLUME = {205},
	PUBLISHER = {American Mathematical Society},
ADDRESS = {Providence, RI},
	YEAR = {2015},
	PAGES = {xvi+343},
	ISBN = {978-1-4704-2024-6},
	MRCLASS = {18D10 (16T20)},
	MRNUMBER = {3242743},
}

@incollection{selinger,
  title = {A survey of graphical languages for monoidal categories},
  author = {Selinger, P.},
  year = 2010,
  booktitle = {New structures for physics},
  location = {Berlin, Heidelberg},
  publisher = {Springer},
  volume = 813,
  pages = {289--355},
  doi = {10.1007/978-3-642-12821-9_4},
  url = {https://www.mscs.dal.ca/~selinger/papers/graphical.pdf},
  urldate = {2018-05-05},
  langid = {english},
  editor = {Coecke, Bob}
}

@book{cwm,
  title = {Categories for the Working Mathematician},
  author = {Mac Lane, S.},
  year = 1998,
  location = {New York, NY},
  publisher = {Springer-Verlag, New York},
  series = {Graduate Texts in Mathematics},
  volume = 5,
  pages = {xii+314},
  doi = {10.1007/978-1-4757-4721-8},
  isbn = {0-387-98403-8},
  url = {http://link.springer.com/10.1007/978-1-4757-4721-8},
  urldate = {2018-05-05},
  edition = {Second},
  mrclass = {18-02},
  mrnumber = 1712872,
  langid = {english},
}

@article{fritz2025categories,
author = {Fritz, Tobias and Lorenzin, Antonio},
title = {Categories of abstract and noncommutative measurable spaces},
journal = {Advances in Mathematics},
volume = {488},
pages = {110793},
year = {2026},
issn = {0001-8708},
doi = {https://doi.org/10.1016/j.aim.2026.110793},
url = {https://www.sciencedirect.com/science/article/pii/S0001870826000150},
}

@inbook{macdonald2003aspects, 
place={Cambridge}, 
series={Encyclopedia of Mathematics and its Applications}, 
title={Aspects of Monads}, 
booktitle={Categorical Foundations: Special Topics in Order, Topology, Algebra, and Sheaf Theory}, 
publisher={Cambridge University Press}, 
author={MacDonald, John and Sobral, Manuela}, 
year={2003}, 
pages={213--268}, 
collection={Encyclopedia of Mathematics and its Applications}
}

@misc {bunke2021lectures,
	author = {Bunke, Ulrich},
	title = {Lecture course on {$C^*$}-algebras and
{$C^*$}-categories},
  year={2021},
	note = {\href{https://bunke.app.uni-regensburg.de/VorlesungCastcats-1.pdf}{bunke.app.uni-regensburg.de/VorlesungCastcats-1.pdf}},
}

@InProceedings{linton1969coequalizers,
author="Linton, F. E. J.",
editor="Eckmann, B.",
title="Coequalizers in categories of algebras",
booktitle="Seminar on Triples and Categorical Homology Theory",
year="1969",
publisher="Springer Berlin Heidelberg",
address="Berlin, Heidelberg",
pages="75--90",
isbn="978-3-540-36091-9"
}

@article{seal2013tensors,
 author = {Seal, Gavin J.},
 title = {Tensors, monads and actions},
 fjournal = {Theory and Applications of Categories},
 journal = {Theory Appl. Categ.},
 issn = {1201-561X},
 volume = {28},
 pages = {403--434},
 year = {2013},
 language = {English},
}

@article{barr2005toposes,
 author = {Barr, Michael and Wells, Charles},
 title = {Toposes, triples and theories},
 fjournal = {Reprints in Theory and Applications of Categories},
 journal = {Repr. Theory Appl. Categ.},
 volume = {2005},
 number = {12},
 pages = {1--288},
 year = {2005},
 language = {English},
}

@article{okamura2015measurement,
    author = {Okamura, Kazuya and Ozawa, Masanao},
    title = {Measurement theory in local quantum physics},
    journal = {Journal of Mathematical Physics},
    volume = {57},
    number = {1},
    pages = {015209},
    year = {2015},
    issn = {0022-2488},
    doi = {10.1063/1.4935407},
    url = {https://doi.org/10.1063/1.4935407},
    eprint = {https://pubs.aip.org/aip/jmp/article-pdf/doi/10.1063/1.4935407/11140738/015209_1_online.pdf},
}

@misc{booth2026quantuminstruments,
  title={Composing Quantum Instruments},
  author={Booth, Robert I. and Leichtle, Dominik and Rice, Alex and Worrall, Kim},
	note = {\href{https://arxiv.org/abs/2606.28291}{arXiv:2606.28291}},
  year={2026}
}

@misc{fritz2026quantuminstrumentmonad,
  title={The quantum instrument monad},
  author={Fritz, Tobias},
	note = {\href{https://arxiv.org/abs/2606.27805}{arXiv:2606.27805}},
  year={2026}
}

@Book{weggeolsen1993Ktheory,
 Author = {Wegge-Olsen, Niels Erik},
 Title = {{{\(K\)}}-theory and {{\(C^*\)}}-algebras: a friendly approach},
 ISBN = {0-19-859694-4},
 Year = {1993},
 Publisher = {Oxford: Oxford University Press},
 Language = {English},
}

@Article{kirchberg1998Cstarfromoperator,
 Author = {Kirchberg, Eberhard and Wassermann, Simon},
 Title = {{{\(C^*\)}}-algebras generated by operator systems},
 FJournal = {Journal of Functional Analysis},
 Journal = {J. Funct. Anal.},
 ISSN = {0022-1236},
 Volume = {155},
 Number = {2},
 Pages = {324--351},
 Year = {1998},
 Language = {English},
 DOI = {10.1006/jfan.1997.3226},
 zbMATH = {1199813},
 Zbl = {0940.46038}
}

@incollection {pelletier1989equationalCstar,
    AUTHOR = {Pelletier, Joan Wick and Rosick{\'y}, Ji{\v r}{\'i}},
     TITLE = {Generating the equational theory of {$C^*$}-algebras and
              related categories},
 BOOKTITLE = {Categorical topology and its relation to analysis, algebra and
              combinatorics ({P}rague, 1988)},
     PAGES = {163--180},
 PUBLISHER = {World Sci. Publ., Teaneck, NJ},
      YEAR = {1989},
      ISBN = {9971-50-966-0},
   MRCLASS = {46L05 (46M15)},
  MRNUMBER = {1047899},
}

@book{adamek1994locallypresentable, 
  place={Cambridge}, 
  series={London Mathematical Society Lecture Note Series}, 
  title={Locally Presentable and Accessible Categories}, 
  publisher={Cambridge University Press}, 
  author={Adamek, J. and Rosicky, J.}, 
  year={1994}, 
  collection={London Mathematical Society Lecture Note Series}
}

@article{davies1970operational,
	author = {Davies, E. Brian and Lewis, John T.},
	date = {1970/09/01},
	doi = {10.1007/BF01647093},
	id = {Davies1970},
	isbn = {1432-0916},
	journal = {Communications in Mathematical Physics},
	number = {3},
	pages = {239--260},
	title = {An operational approach to quantum probability},
	url = {https://doi.org/10.1007/BF01647093},
	volume = {17},
	year = {1970}}

@article{ozawa1984quantummeasuring,
    author = {Ozawa, Masanao},
    title = {Quantum measuring processes of continuous observables},
    journal = {Journal of Mathematical Physics},
    volume = {25},
    number = {1},
    pages = {79-87},
    year = {1984},
    month = {01},
    issn = {0022-2488},
    doi = {10.1063/1.526000},
    url = {https://doi.org/10.1063/1.526000},
    eprint = {https://pubs.aip.org/aip/jmp/article-pdf/25/1/79/19265867/79_1_online.pdf},
}

@article{khoshkam2002categorical, 
title={Categorical constructions in {C}*-algebra theory}, 
volume={73}, 
DOI={10.1017/S1446788700008491}, 
number={1}, 
journal={Journal of the Australian Mathematical Society}, 
author={Khoshkam, M. and Tavakoli, J.}, 
year={2002}, 
pages={97--114}
}

@article{shirazi2026commutativity,
 author = {Shirazi, Zev},
 title = {Commutativity and {Kleisli} laws for codensity monads of probability measures},
 fjournal = {Theory and Applications of Categories},
 journal = {Theory Appl. Categ.},
 issn = {1201-561X},
 volume = {45},
 pages = {461--500},
 year = {2026},
 language = {English},
 url = {www.tac.mta.ca/tac/volumes/45/14/45-14abs.html},
}

@book {rudin1987analysis,
    AUTHOR = {Rudin, Walter},
     TITLE = {Real and complex analysis},
   EDITION = {Third},
 PUBLISHER = {McGraw-Hill Book Co., New York},
      YEAR = {1987},
     PAGES = {xiv+416},
      ISBN = {0-07-054234-1},
   MRCLASS = {00A05 (26-01 30-01 46-01)},
  MRNUMBER = {924157},
}
